\documentclass[11pt]{amsart}
\usepackage{amsmath,amsfonts,amssymb}
\usepackage{enumerate,setspace,graphics,graphicx}
\usepackage{latexsym,amscd,subfigure,comment}
\usepackage[margin=1in]{geometry}
\usepackage{fullpage}
\usepackage[usenames]{xcolor}
\usepackage{verbatim}

\usepackage{cancel}
\usepackage{bm} 
\usepackage{hyperref}
\usepackage{yfonts}

\newtheorem{thm}{Theorem}[section]
\newtheorem{cor}[thm]{Corollary}
\newtheorem{lem}[thm]{Lemma}
\newtheorem{lemma}[thm]{Lemma}
\newtheorem{prop}[thm]{Proposition}
\newtheorem*{thmA}{Theorem A}
\newtheorem*{thmB}{Theorem B}
\newtheorem*{thmC}{Theorem C}
\newtheorem*{thmD}{Theorem D}

\theoremstyle{definition}
\newtheorem{defn}[thm]{Definition}

\newtheorem{remark}[thm]{Remark}
\newtheorem{ex}[thm]{Example}
\newtheorem{nota}[thm]{Notation}
\newtheorem{notarem}[thm]{Notation and Remarks}

\definecolor{dgreen}{RGB}{0,128,0}

\newcommand{\R}{\mathbf R}
\newcommand{\Z}{\mathbf Z}

\newcommand{\mix}{\operatorname{mix}}
\newcommand{\Stek}{\operatorname{Stek}}

\newcommand{\Om}{\Omega}

\newcommand{\pa}{\partial}

\newcommand{\Lc}{\mathcal{L}}

\newcommand{\A}{\mathcal A}

\newcommand{\B}{\mathbf B}
\newtheorem{ques}[thm]{Question}

\newcommand{\restr}[1]{\lower0.4ex\hbox{$|$}\lower0.7ex
	\hbox{$\scriptstyle{#1}$}}
\definecolor{olive}{RGB}{128,128,0}

\definecolor{orange}{rgb}{1,0.5,0}

\definecolor{darkpastelpurple}{rgb}{0.59,0.44,0.84}

\newcommand{\la}{\langle}

\newcommand{\ra}{\rangle}

\newcommand{\psom}{\partial_S\Om}
\newcommand{\pnom}{\partial_N\Om}

\newcommand{\pst}{\partial_*}

\newcommand{\disk}{{\mathbb{D}}}
\newcommand{\hD}{{\mathbb{D}}_{\scriptscriptstyle{\frac{1}{2}}}}
\newcommand{\qD}{{\mathbb{D}}_{\scriptscriptstyle{\frac{1}{4}}}}
\newcommand{\bell}{\bar{\ell}}
\newcommand{\pdom}{\partial_D\Om}

\newcommand{\om}{\Omega}

\newcommand{\mL}{\mathbf L}

\begin{document}

\title[Steklov and mixed Steklov problems on surfaces]{Applications of possibly hidden symmetry to Steklov and mixed Steklov problems on surfaces}

\author[Arias-Marco]{Teresa Arias-Marco}
\address{Department of Mathematics, University of Extremadura,  Av. de Elvas s/n, 06006 Badajoz, Spain}
\email{ariasmarco@unex.es}
\author[Dryden]{Emily B. Dryden}
\address{Department of Mathematics, Bucknell University, Lewisburg, PA 17837, USA}
\email{emily.dryden@bucknell.edu}
\author[Gordon]{Carolyn S. Gordon}
\address{Department of Mathematics, 6188 Kemeny, Dartmouth College, Hanover, NH  03755, USA}
\email{csgordon@dartmouth.edu}
\author[Hassannezhad]{Asma Hassannezhad}
\address{University of Bristol,
School of Mathematics,
Fry Building,
Woodland Road,
Bristol, 
BS8 1UG, U.K.}
\email{asma.hassannezhad@bristol.ac.uk}
\author[Ray]{Allie Ray}
\address{Birmingham-Southern College, Department of Mathematics, 900 Arkadelphia Rd, Birminham, AL 35254, USA}
\email{adray@bsc.edu}
\author[Stanhope]{Elizabeth Stanhope}
\address{Department of Mathematical Sciences, Lewis \& Clark College, Portland, OR 97219, USA}
\email{stanhope@lclark.edu}


\date{\today}

\begin{abstract}
We consider three different questions related to the Steklov and mixed Steklov problems on surfaces.  These questions are connected by the techniques that we use to study them, which exploit symmetry in various ways even though the surfaces we study do not necessarily have inherent symmetry.

In the spirit of the celebrated Hersch-Payne-Schiffer and Weinstock inequalities for Steklov eigenvalues, we obtain a sharp isoperimetric inequality for the mixed Steklov eigenvalues considering the interplay between the eigenvalues of the mixed Steklov-Neumann and Steklov-Dirichlet eigenvalues. 

In 1980, Bandle showed that the unit disk maximizes the $k$th nonzero normalized Steklov eigenvalue on simply connected domains with rotational symmetry of order $p$ when $k\le p-1$. We discuss whether the disk remains the maximizer in the class of simply connected rotationally symmetric domains when $k\geq p$. In particular, we show that for $k$ large enough, the upper bound converges to the Hersch-Payne-Schiffer upper bound.

We give full asymptotics for mixed Steklov problems on arbitrary surfaces, assuming some conditions at the meeting points of the Steklov boundary with the Dirichlet or Neumann boundary. 
\end{abstract}

\maketitle

\section{Introduction}

The Steklov eigenvalue problem was introduced at the turn of the 20th century by Vladimir Steklov \cite{S02}.  The relationship between the Steklov eigenvalues and the geometry of the domain under consideration has been studied for over fifty years (e.g., \cite{ADGHRS,BKPS,GP17,HL20,LPPS,W}).  For a physical interpretation of these eigenvalues, one may think of a planar domain as representing a vibrating free membrane with some mass density concentrated along the boundary. The Steklov eigenvalues are the squares of the natural vibration frequencies of this membrane. 

Throughout the paper, unless otherwise stated, \label{omega}$(\Om,g)$ will denote an orientable, connected, compact Riemannian surface with nonempty Lipschitz boundary. By Lipschitz boundary, we mean that $(\Om,g)$ is isometric to a domain with Lipschitz boundary in a complete Riemannian manifold.  The Steklov problem on $(\Om, g)$ is given by 
\begin{equation}
\begin{cases}
 \Delta_g f = 0 & {\rm in} \ {\Om}, \\
      \partial_\nu  f = \sigma  f & {\rm on} \ \partial  {\Om},
\end{cases}
\end{equation}
where $\Delta$ is the Laplace-Beltrami operator acting on functions on $\Om$ and $\nu$ is the unit outward normal vector field along the boundary (defined almost everywhere).  
The  Steklov spectrum $\Stek(\Om)$ consists of a discrete sequence of nonnegative real numbers 
\[0=\sigma_0(\Omega)< \sigma_1(\Omega)\le \sigma_2(\Omega)\le \cdots \nearrow +\infty\]
and $\infty$ is the only accumulation point. 

In addition to the Steklov problem, we also consider mixed eigenvalue problems. Decompose $\pa\Om$ into three parts, $\pa\Om=\pa_S\Om\sqcup\pa_N\Om \sqcup \pa_D\Om$, where each part is a finite union of topological segments and/or circles and $\pa_S\Om$ is required to be nonempty. The general mixed eigenvalue problem is given by
\begin{equation}\label{genmixdef}
\begin{cases}
 \Delta_g f = 0 & {\rm in} \ {\Om}, \\
      \partial_\nu f = \sigma  f & {\rm on} \ \partial_S  {\Om},\\
       \partial_\nu f = 0& {\rm on} \ \partial_N  {\Om},\\
         f = 0& {\rm on} \ \partial_D  {\Om}.
\end{cases} 
\end{equation}
 When $\pa_D\Om$ is empty, we obtain the Steklov-Neumann problem with spectrum $\Stek_N(\Omega)$ written as
\[0=\sigma_0^N(\Om)<\sigma_1^N(\Om)\le\cdots \nearrow +\infty.\]
Similarly when $\pa_N \Om$ is empty, we obtain the Steklov-Dirichlet problem with spectrum $\Stek_D(\Omega)$ written as 
 \[0<\sigma_0^D(\Om)\le\sigma_1^D(\Om)\le \cdots \nearrow +\infty.\]
Note that the indexing of the Steklov-Dirichlet eigenvalues starts with 0. This might be unconventional, but it helps avoid possible confusion when we compare $\Stek_N(\Om)$ and $\Stek_D(\Om)$.  When we wish to consider both the Steklov-Neumann and Steklov-Dirichlet eigenvalue problems at once, it will be convenient to let $\pa_*\Om$ denote either $\pa_D\Om$ or  $\pa_N\Om$, as needed. 

The Steklov-Neumann problem is also known as the sloshing problem for domains with certain geometric characteristics.  The sloshing problem is essential in hydrodynamics and can be used to study thermodynamics where the ``walls'' $\partial_N {\Om}$ of the domain are insulated (see \cite{LPPS}). If the temperature of the walls is kept at zero, the thermodynamics problem becomes a particular case of the Steklov-Dirichlet problem (see \cite{HP68}). \\

 The purpose of this paper is threefold:
 
 \begin{itemize}
 \item Fixing a boundary decomposition $\pa\om =\psom \cup \pa_*\om$ on a surface of genus zero, we prove a sharp isoperimetric inequality for $\min\{\sigma_k^N(\om,g), \sigma_{k-1}^D(\om,g)\}L(\pa_S\Om)$ in the spirit of  the celebrated Hersch-Payne-Schiffer inequality \eqref{HPS1}.   
 \item Motivated by Weinstock's observation that the Euclidean disk maximizes $\sigma_1(\Om)L(\pa_S\Om)$ among planar domains, we examine the extent to which it continues to be a maximizer for higher eigenvalues over families of domains with symmetry. Our results build off of an isoperimetric inequality of Bandle \cite{bandle}. 
 \item We give the full asymptotics for mixed Steklov eigenvalue bounds on surfaces of arbitrary genus and number of boundary components under constraints at all meeting points of $\pa_S\Om$ with $\pa_N\Om$ and $\pa_D(\om)$.  We give several applications, e.g., to the inverse spectral problem. 
 \end{itemize}

Before describing our main results in detail, we make a few observations connecting various forms of the Steklov problem. The Steklov problem and  Steklov-Neumann problem can be viewed as special cases of the weighted Steklov problem:  
 \begin{equation}\label{first def}
\begin{cases}
 \Delta_g f = 0 & {\rm in} \ {\Om}, \\
      \partial_\nu  f = \sigma \rho f & {\rm on} \ \partial  {\Om},
\end{cases}
\end{equation}
 where  $\rho\in L^\infty(\pa\Om)$ is a non-negative weight function. 
When the weight function $\rho\in L^\infty(\pa\Om)$ is identically one on $\pa \Om$, \eqref{first def} is the Steklov eigenvalue problem.  Suppose for the moment that $\rho$ takes on only the values 0 and 1. Denoting the support of $\rho$ by $\pa_S\Om$, and the set of boundary points at which $\rho$ is zero by  $\pa_N\Om\ne\emptyset$  gives the Steklov-Neumann problem. 

Using recent results of Dorin Bucur and Micka\"el Nahon \cite{BN20}, we show in Subsection \ref{smooth to lipshitz} that any Steklov eigenvalue bound that is valid for surfaces with smooth boundary is also valid for surfaces with Lipschitz boundary and also for weighted Steklov eigenvalue problems on surfaces. Mikhail Karpukhin and Jean Lagac\'e \cite{KL22} independently obtained the same result for manifolds of all dimensions. Thus for the simplicity of the presentation, we state our results as well as previous results (including those that originally assumed stronger regularity on $\pa\Om$) for unweighted Steklov problems on surfaces with Lipschitz boundary, keeping in mind that the statements still hold if one  considers the weighted Steklov problem with a non-negative, nontrivial  $L^\infty$ weight.  Only the results on the asymptotics of eigenvalues require the  weight  to be smooth and strictly positive on $\pa_S\Om$, and also require $\psom$ to be smooth.

\subsection{Bounds for mixed eigenvalue problems} The literature contains many isoperimetric bounds for normalized Steklov eigenvalues on surfaces.   By the remarks above, all such bounds are also valid for mixed Steklov-Neumann problems. 

Robert Weinstock~\cite{W} (see also \cite[Thm.~1.3]{GP09}) initiated the study of upper bounds for Steklov eigenvalues of surfaces in 1954 by deriving what is now called the Weinstock inequality. For $\Om$ a simply-connected Riemannian surface with boundary, one has
 \begin{equation}\label{weinstock}
 \sigma_1(\Om) L(\partial \Om)\leq 2\pi
 \end{equation}
with equality if and only if $\Om$ is a Euclidean disk. Joseph Hersch, Lawrence Payne, and Menahem Schiffer \cite{HPS} generalized the Weinstock inequality to higher eigenvalues of simply-connected surfaces:
 \begin{equation}\label{HPS1}
 \sigma_k(\Om)L(\partial \Om)\leq 2\pi k.    
 \end{equation}
 Both these inequalities were originally stated only in the setting of plane domains and have been generalized and improved in recent decades, e.g. \cite{Kar17,YY17,GP,GP2012}.
 
 For $k\geq 2$, Alexandre Girouard and Iosif Polterovich \cite{GP} proved that the Hersch-Payne-Schiffer inequality \ref{HPS1} is sharp, and Ailana Fraser and Richard Schoen \cite{FS2020} proved that it is strict at least among surfaces with smooth boundary.

For surfaces of genus zero with an arbitrary number of boundary components, more recent work of \cite{GKL, Kok14} yields the strict inequality
\begin{equation}\label{HPS2}
 \sigma_k(\Om)L(\partial \Om)<8\pi k.    
 \end{equation}

When $\om$ is a topological annulus, Fraser and Schoen \cite{FS2} showed that the intersection of the so-called critical catenoid and the Euclidean ball $\B^3$ maximizes  $\sigma_1(\Om)L(\pa\Om)$.  They then obtained the improved upper bound of approximately $\frac{4\pi}{1.2}$ for $\sigma_1(\om,g)L_g(\pa\om)$ for Riemannian metrics on the annulus.
 
As noted above, Inequalities~\eqref{HPS1}, respectively \eqref{HPS2}, immediately yield the same bounds for normalized mixed Steklov-Neumann eigenvalue bounds on simply-connected Riemannian surfaces, respectively Riemannian surfaces of genus zero.  One also knows from the variational characterization of eigenvalues that if one fixes the Steklov part of the boundary and imposes either Neumann or Dirichlet conditions on the remaining part, then $\sigma_k^N(\om)\leq \sigma_k^D(\om)$.   
In the spirit of the Weinstock and Hersch-Payne-Schiffer inequalities, we prove:

\begin{thmA} 
Assume that $\Omega$ has genus 0. Let $\partial \Om= \psom \sqcup \pa_*\Om$ be a nontrivial decomposition of $\pa\Om$. We assume that all the connected components of $\pa_*\Om$ are  contained in  a single connected component of $\pa\Om$.  Then
\begin{equation}\label{thma}\min\{\sigma_k^N(\Om), \sigma_{k-1}^D(\Om)\}L(\psom)< 4(2k-1)\pi, \qquad  k\ge1.\end{equation}
Moreover, if equality holds, then $\sigma_{k-1}^D(\Om)=\sigma_k^N(\Om)$. 
 Inequality\eqref{inq:gkl} is sharp, but the number of boundary components of any sequence of domains approaching the upper bound cannot remain bounded.  
 
If $\Omega$ is simply connected and if $\psom$ and $\pa_*\om$ are connected,  inequality \eqref{thma}  can be improved to the sharp inequality:
 \begin{equation}\label{thma2}\min\{\sigma_k^N(\Om), \sigma_{k-1}^D(\Om)\}L(\psom)\leq (2k-1)\pi,\qquad  k\ge1.\end{equation}
When $k=1$, equality in \eqref{thma2} is achieved for the flat half disk discussed in Example~\ref{models}.
\end{thmA}
We refer to Theorem \ref{general minpi} for a more complete statement of  Theorem A. 

When $k=1$ and either (i) $\om$ is a topological annulus and each of $\psom$ and $\pa_*\om$ is one of the boundary circles or (ii) $\om$ is a topological disk and each of $\psom$ and $\pa_*(\om)$ consists of two intervals, we  find a maximizing metric $g_0$ for $\min\{\sigma_k^N(\Om), \sigma_{k-1}^D(\Om)\}L(\psom)$.   In the first case, the maximum value is realized by a portion of the critical catenoid bounded by the equatorial plane and $\mathbf{B}^3$, and in the second case by a portion of the critical catenoid bounded by $\mathbf{B}^3$ and a plane through the catenoid's axis.
In both cases, the part of the boundary on which we impose Neumann or Dirichlet conditions consists of the intersection of the critical catenoid with the cutting plane.

We also discuss some settings in which we can make stronger conclusions. For example, if the domain $\Om$ satisfies the \textit{weak John's condition} (see Definition~\ref{def.john}), then Theorem A along with a result of \cite{BKPS} imply (Proposition \ref{prop.john})
$\sigma_k^N(\Om)L(\pa\Om) \le (2k-1)\pi.$
We show that this inequality is sharp when $k=1$.\\ 

To prove Theorem A, we first establish it under the additional hypothesis that $\om$ is isometric to a quotient of a surface $\Sigma$ by a reflection symmetry; we refer to this additional hypothesis as the ``Lipschitz doubling condition''.   The Steklov-Neumann (respectively, Steklov-Dirichlet) eigenvalues of $\Om$ correspond to Steklov eigenvalues on $\Sigma$ whose eigenfunctions are invariant (respectively, anti-invariant) under the reflection. For the general case, we approximate $\om$ by surfaces $\om_j$ for which $\pa_*\om$ is polygonal; each of these surfaces is isospectral to a surface $\om_j'$ satisfying the Lipschitz doubling condition and thus the theorem holds for the $\om_j$'s.  A stability result of Bucur and Nahon \cite{BN20} and domain monotonicity of eigenvalues are used to complete the argument in the general case.

\subsection{Upper bounds for domains with symmetry} 

Although Weinstock showed that $\sigma_1(\Om)L(\pa \Om)$ is maximized by a Euclidean disk $\disk$ among all simply connected planar domains, the disk fails to maximize higher eigenvalues.   However, Catherine Bandle showed that if one restricts attention to simply connected domains with rotational symmetry of order $p$, then the Euclidean disk continues to maximize $\sigma_k(\Om)L(\pa \Om)$ for all $k\leq p-1$.  We pursue this line of inquiry, asking whether Bandle's Theorem gives the maximum $k$ for which the Euclidean disk maximizes  $\sigma_k(\Om)L(\pa \Om)$.

Let  $\mathcal{D}_p$ be the class of all 
simply connected planar domains  with Lipschitz  boundary and rotational symmetry of order $p\ge2$ centered at the origin and set
$$\Sigma_k^{p}=\sup_{\Om\in \mathcal{D}_p}\sigma_k(\Om)L(\pa \Om).$$ 
Let $K(p)$ denote \[K(p):=\{k\in \Z^+: \sigma_k(\disk )L(\partial \disk )=\Sigma_k^{p}\}\] 
Because the domains constructed by Girouard and Polterovich's \cite{GP} in their proof of sharpness for the Hersch-Payne-Schiffer inequality lie in $\mathcal{D}_2$, they can be used to conclude that $K(2)\cap[2,\infty)=\emptyset$.   Motivated by their construction, we show that
\begin{thmB} $\min  K(p)=p$ when $p$ is even and $\min  K(p)\le2p$ when $p$ is odd. 
 \end{thmB}
For a finite number of odd integers $p$, we apply numerical computations of Eldar Akhmetgaliyev, Chiu-Yen Kao and Braxton Osting \cite{AKO} (see also \cite{Bo17}) to obtain the stronger result that $\min  K(p)=p$.   
Theorem B is an immediate consequence of Theorem \ref{bndlesharp} and Corollary \ref{BandleLimit}, which give more descriptive information about the set $K(p)$ and yield the following asymptotic result:
\begin{thmC}
For every $p$, we have 
$\displaystyle{\lim_{k\to\infty}\,\frac{\Sigma^p_k}{2\pi k}=1}$.
\end{thmC}

\subsection{Eigenvalue asymptotics for mixed Steklov problems} 

We now shift from upper bounds on eigenvalues to eigenvalue asymptotics for mixed Steklov problems.  In the Steklov-Neumann setting, such asymptotics have a long history going back to \cite{San55}.   Michael Levitin, Leonid Parnovski, Iosif Polterovich, and David Sher \cite{LPPS} obtained two term-asymptotics for the Steklov-Neumann and Steklov-Dirichlet problems when $\Om$ is a simply connected planar domain with Lipschitz boundary and boundary decomposition $\pa\om=\pa_S\Om \cup \pa_*\om$, where $\psom$ is a straight line segment and where $\psom$ and $\pa_*\om$ meet at angles $ \alpha ,\beta\in(0, \frac{\pi}{2}]$. Verifying a conjecture by David William Fox and James R. Kuttler \cite{FK}, they show 
\begin{equation}\label{LPPS_asymptotics}
\sigma_k^{N}(\Om) L(\pa_S \Om) = \pi \left(k + \frac{1}{2}\right) - \frac{\pi^2}{8}\left(\frac{1}{\alpha} + \frac{1}{\beta}\right) + r(k),\qquad k \rightarrow \infty
\end{equation}
where   $r(k)$ is a function  of order at least $o(1)$. (We have adjusted their formula to our indexing conventions.) A similar two-term asymptotic holds for the Steklov-Dirichlet eigenvalue. The only difference is that the second term has a positive sign.  The order of $r(k)$ depends on the meeting angles of $\pa_S\Om$ and $\pa_*\Om$  and the geometry of the boundary near the meeting points. 
The decay rate of $r(k)$ can be improved to $r(k)=o(e^{-ck})$, for some constant $c>0$,  when $\alpha=\beta=\frac{\pi}{2}$  
(see \cite[Propositions 1.3 and 1.8]{LPPS}.
 For more details  see \cite{LPPS,LPPS2} and references therein for earlier work.
 
Now let $\Om$ be a Riemannian surface with Lipschitz boundary of arbitrary genus and arbitrary number of boundary components. We consider the mixed Steklov-Neumann-Dirichlet problem where $\pa\Om=\pa_S\Om\sqcup\pa_N\om\sqcup\pa\Om_D$ gives a  decomposition into three parts.  (We allow the case that $\pnom$ or $\pdom$ is empty.)  Steklov, Neumann, and Dirichlet boundary conditions are imposed on   $\pa_S\Om$, $\pa_N\om$, and $\pa\Om_D$ respectively. 
 
\begin{thmD} Let $\Om$ be a Riemannian surface with Lipschitz boundary and write $\pa\Om=\pa_S\Om\sqcup\pa_N\om\sqcup\pa\Om_D$ as above.     In addition, assume $\psom$ is smooth and 
at any common endpoint $p$ of $\psom$ with either of $\pnom$ or $\pdom$, a small segment of $\pnom$, respectively $\pdom$, emanating from $p$ is a geodesic segment orthogonal to $\psom$ at $p$. 
 Then the full asymptotics of the mixed Steklov-Dirichlet-Neumann problem is given by the corresponding mixed Steklov-Dirichlet-Neumann asymptotics of a disjoint union of disks, half disks, and quarter disks with perimeter equal to the length of one of the boundary components of $\pa_S\Om$.  
\end{thmD}
For the precise statement of Theorem D including the boundary conditions that are imposed on the disks, half disks, and quarter disks, we refer to Theorem \ref{asymptotic-dir-Neu}. By ``full asymptotics'' we mean that the error is of order $O(k^{-\infty})$ where $k$ is the index of the eigenvalues. The assumption that small segments of $\pa_*\Om$ near any common endpoints with $\pa_S\Om$ are geodesic cannot be completely removed. By a result of Anthony M. J. Davis \cite{Davis} (see also \cite{Ursell,Davis2}), the curvature of $\pa_*\Om$ near the endpoints appears in the second term of the asymptotics. 

Theorem D can be viewed as an extension of the result of Girouard, Parnovski, Polterovich and Sher \cite{GPPS.14} on the full asymptotics of the   Steklov spectrum for Riemannian surfaces  with smooth boundary, and the result of the authors \cite{ADGHRS} on the full asymptotics of the Steklov spectrum for orbisurfaces.

Theorem D, more precisely Theorem \ref{asymptotic-dir-Neu}, has a number  of interesting consequences for Steklov-Neumann or Steklov-Dirichlet problems of which we list a few below. More details and further consequences appear in Section \ref{sec:4}.  
\begin{enumerate}
    \item (Inverse spectral result.) Consider any mixed Steklov-Neumann or Steklov-Dirichlet problem satisfying the hypotheses of Theorem D.  Then the spectrum determines the number of components of $\pa_S\Om$ that are topological circles $m$ and the number $n$ that are topological intervals, and it determines their lengths up to an equivalence relation.  See Corollary \ref{cor: inverse} and Remark \ref{rem: 5.7} for details.
    \item In the notation of item (1), there exists $k_\Om>0$ such that the multiplicity of $\sigma_k^N(\Om)$ and $\sigma_k^D(\Om)$ is at most $2m+n$ for every $k\ge k_\Om$.
    \item Considering the Steklov-Neumann and Steklov-Dirichlet problems on the same domain $\om$ with $\psom$ fixed, we have 
 \[ \sigma_{k+m}^N(\Omega)- \sigma_k^D(\Omega)=O(k^{-\infty})\]
 where $m$ is the number of components of $\psom$ that are topological intervals. 
\end{enumerate}

The proof of Theorem D follows a sequence of steps similar to those in Theorem A, although the methods used in each step are different. We first address the case in which the surface satisfies the ``Lipschitz doubling condition'' applying results of \cite{ADGHRS}.  The case in which each component of $\pnom\sqcup\pdom$ is polygonal then follows.  The latter case along with domain monotonicity of eigenvalues are used to address the general case.

The paper is organized as follows. In Section \ref{chap.background_rev}, we cover the necessary background and preliminaries. This includes, in particular, the background on extending eigenvalue bounds from surfaces with smooth boundary to ones with Lipschitz boundaries and also introduces the Lipschitz doubling condition referred to above.  
In Section \ref{bounds_mixed}, we prove Theorem A and its consequences.  Section \ref{Bandle} investigates the sharpness of Bandle's Theorem for domains with symmetry.
 The last section deals with the full eigenvalue asymptotics for the mixed Steklov-Neumann-Dirichlet eigenvalues and  its corollaries.

\subsection*{Acknowledgements}
We would like to thank B. Bogosel, M. Karpukhin, O. Klurman,  and J. Lagac\'e for helpful discussions during this project. We are also grateful for our time spent in the Research in Groups programme at the International Centre for Mathematical Sciences in Edinburgh, UK, for our visit to the Department of Mathematics of Bucknell University in Lewisburg, USA supported by a grant from the Simons Foundation (210445 to Emily B. Dryden), as well as for our stay at the Department of Mathematics of the University of Extremadura in Badajoz, Spain under the GR18001 grant funded by Junta de Extremadura and Fondo Europeo de Desarrollo Regional.  T. A.  is partially supported by grant PID2019-10519GA-C22 funded by AEI/10.13039/501100011033 and by grant GR21055 funded by Junta de Extremadura and Fondo Europeo de Desarrollo Regional, and A. H. is partially supported by EPSRC grant EP/T030577/1.

\section{Preliminaries}\label{chap.background_rev}

\subsection{Background}\label{subsec.back}

We will consider Steklov and mixed problems on $(\Om,g)$, possibly with a non-negative weight $\rho$ on $\pa \Om$.   Since we assume that $\Om$ is connected, zero occurs in the spectrum with multiplicity one except in the case of mixed problems in which Dirichlet conditions are imposed on part of the boundary.   
\begin{nota}\label{nota:basic} When considering mixed problems on $(\Omega,g)$, we will write the boundary decomposition as $\partial\Om=\psom \sqcup \pnom \sqcup \pdom$
where Steklov, Neumann, and Dirichlet boundary conditions are placed on $\psom$, $\pnom$, and $\pdom$, respectively.    (Either of $\pnom$ or $\pdom$ may be empty.) Our standing assumption is that 
each of $\psom$, $\pnom$, and $\pdom$ is a  finite disjoint union of intervals or circles and any two of $\psom$, $\pnom$ and $\pdom$ intersect only at endpoints of the intervals.  When we want to consider mixed Steklov-Neumann and mixed Steklov-Dirichlet problems simultaneously, we will often write $\partial\om=\psom \sqcup \partial_*\om$ where either Neumann or Dirichlet conditions will be placed on $\partial_*\om$. 
In the case of mixed problems, we will always assume, even if not stated, that $\psom$ and at least one of $\pnom$ and $\pdom$ are non-trivial. 

We will use the notation $\Stek_N(\om)$ and $\Stek_D(\om)$ for the spectra of mixed Steklov-Neumann and Steklov-Dirichlet problems, respectively, and the notation $\Stek_{DN}$ for fully mixed Steklov-Dirichlet-Neumann problems.  We will also write $\Stek_{mix}(\om)$ when we want to consider all three types of problems simultaneously (see \S \ref{sec:4}).

When referring to weighted Steklov or mixed Steklov problems on $(\Om,g)$, weight $\rho\in L^\infty(\partial_S\Om)$ is always assumed to be non-negative and not identically zero (i.e. to be positive on a subset of positive measure). We will say a weight $\rho$ is trivial if $\rho\equiv 1$.
\end{nota}

We first recall the well-known variational characterizations of eigenvalues for the weighted Steklov problem and mixed problems with (possibly trivial) weight $\rho$ and conformal invariance of the eigenvalues in dimension two. 
Let $(\Om,g)$ have boundary decomposition $\pa\Om=\psom\sqcup\pa_*\Om$.  
\begin{enumerate}
\item  For the Steklov problem (when $\pa_S\Om=\pa\Om$) or the Steklov-Neumann problem with possibly trivial weight, we have
\begin{equation}\label{VarChar}
\sigma_k(\Om,g,\rho)=\inf_{V_{k+1}}\sup_{0\ne f\in V_{k+1}}\frac{\int_\Om|\nabla_g f|^2 dA_g}{\int_{\pa_S\Om}f^2 \rho ds_g},
\end{equation} 
where $V_{k+1}$ varies over all $k+1$-dimensional subspaces of $C^\infty(\Om)$.

\item For mixed Steklov-Dirichlet or Steklov-Neumann-Dirichlet problems with possibly trivial weight $\rho$, the $k$th eigenvalue is given by
\begin{equation}\label{VarChar2}
\inf_{E_{k+1}}\sup_{0\ne f\in E_{k+1}}\frac{\int_\Om|\nabla_g f|^2 dA_g}{\int_{\pa_S\Om}f^2 \rho ds_g},
\end{equation} 
where $E_{k+1}$ consists of all $k+1$-dimensional subspaces of $C^\infty(\Om)$ consisting of functions that are supported on the complement of  $\pdom$ (or its proper subset) where the Dirichlet boundary condition is imposed.
\item For all of the problems above, the eigenvalues are invariant under conformal changes of metric on $\Om$, provided that the conformal factor is identically one on $\pa_S\Om\cap {\rm supp}(\rho)$.    Moreover, the numerator in the variational quotient (the Dirichlet energy) is invariant under all conformal changes of the metric.
\end{enumerate}

\begin{remark}\label{item:smoothpos}\label{rem:linking_the_probs}
Let $(\Om,g)$ be an orientable, connected, compact Riemannian surface with nonempty smooth boundary and suppose that the weight function $\rho$ on $\pa\Om$ is strictly positive and smooth. Then the weighted Steklov eigenvalue problem on $(\Om, g)$ is isospectral to an unweighted Steklov problem on $(\Om,\tilde g)$ for some metric $\tilde g$ in the conformal class of $g$. That is, $\sigma_k(\Om,g, \rho)=\sigma_k(\Om,\tilde g, 1)$ for every $k \ge0$.  
Indeed, we can choose a metric $\tilde g$ conformal to $g$ on $\Om$ such that $\tilde g=\rho^2 g$ on $\pa_S\Om$. Then $ds_{\tilde g}=\rho ds_g$ and the statement follows from the min-max principle~\eqref{VarChar} and the conformal invariance of the Dirichlet energy.

A similar statement holds for mixed problems in which $\pa_S\Om$ is smooth and $\rho$ is strictly positive on $\pa_S\Om$, including at any endpoints of the arcs making up $\pa_S\Om$.   
\end{remark}

\subsection{Lipschitz boundary and regularization}\label{sect.lipschitz}

 We first review the notion of domains with Lipschitz boundaries (see \cite[Appendix A]{MT99}, see also \cite[Definition 2.4.5]{HP05}).  Subsequently, we discuss the approximation of Lipschitz domains by smooth domains and recall a recent result of Bucur and Nahon (see \cite{BN20}) addressing convergence of Steklov eigenvalues under such approximations.

\begin{defn}\label{def.lip}
\begin{enumerate}[1.]
\item Let $(M,g)$ be a Riemannian manifold of dimension $n$ and let $\Om$ be a bounded domain in $M$. 
\begin{enumerate}[a.]
\item A chart $(U_p,h_p)$ on $M$ centered at a point $p\in \pa\Om$ (i.e., $h_p(p)=0$) will be called a \emph{special chart} if it satisfies the following three properties:
\begin{enumerate}
\item $h_p$ is bi-Lipschitz relative to the Riemannian structure on $U_p$ and the Euclidean structure on $h_p(U_p)$;
\item There exists an interval $I_p$ centered about $0$ and a positive real number $r_p$ such 
that $h_p(U_p)=B_{r_p}^{n-1}(0)\times I_p$ where $B_{r_p}^{n-1}(0)$ is the open ball of radius $r_p$ in $\R^{n-1}$ centered at the origin;
\item $h_p(U_p\cap \pa \Om)$ is the graph of a Lipschitz function $\phi_p: B_{r_p}^{n-1}(0)\to I_p$ with $\phi_p(0)=0$ and $$h_p(U_p\cap \Om)=\{(x,t)\in B_{r_p}^{n-1}(0)\times I_p: t>\phi_p(x)\}.$$
\end{enumerate}
\item A domain $\Om\subset M$ is said to have a \emph{Lipschitz boundary} if for every $p\in \pa\Om$, there exists a special chart $(U_p,h_p)$ in $M$ centered at $p$.
\item We say $\Om\subset M$ is a domain with \emph{uniformly Lipschitz boundary} if the collection of special charts $\{(U_p,h_p)\}_{p\in\pa\Om}$ can be chosen so that the following conditions are satisfied:  the bi-Lipschitz constants of the maps $h_p$ and the Lipschitz constants of the functions $\phi_p$ are uniformly bounded independently of $p$,  and $r_p$ and the length of $I_p$ are bounded below by positive constants independent of $p$.
\end{enumerate}
\item\label{stand}
We will say that $(\Om,g)$ is a \emph{Riemannian manifold with Lipschitz boundary} if it is isometric to a domain with Lipschitz boundary in a complete Riemannian manifold.  Note that every Riemannian manifold with smooth boundary satisfies this condition. 
\end{enumerate}
\end{defn}
\begin{remark}\label{rem.weak}~ 
\begin{enumerate}

\item  Definition \ref{def.lip}(2) is a priori weaker than the standard definition of a domain with Lipschitz boundary (see e.g., \cite[Definition 2.4.5]{HP05}) in the Euclidean setting.   In that setting, the chart maps $h_p:U_p \to B_{r_p}^{n-1}(0)\times I_p$ are normally assumed to be Euclidean isometries.  However, by \cite[Definition 2.5 and Theorem 4.1]{HMT07}, Definition \ref{def.lip} is in fact equivalent to the usual definition in the Euclidean setting.
\item  Let $\Om$ be a domain with Lipschitz boundary in an $n$-dimensional Riemannian manifold and let $p\in \pa\Om$.  Given any chart $(V,h)$ centered at $p$, then by \cite[Theorem 4.3]{HMT07}, there exists an open set $U\subset V$ containing $p$ and a Euclidean isometry $\tau$ of $\R^n$ such that $(U, \tau\circ h)$ is a special chart.   
\item \label{part3}
We will only consider domains $\Om$ with Lipschitz boundary that have compact closure.  Such domains always have uniformly Lipschitz boundary.    In fact, starting with any finite collection $\mathcal{C}$ of special charts whose domains cover $\pa\Om$, there exists $r_0>0$ and an open interval $I_0$ about $0$  in $\R$ satisfying the following condition:   For every $p\in \pa\Om$, there exists $(U,h)\in\mathcal{C}$ such that $p\in U$ and the image of the map $h_p:=h-h(p)$ contains the cylinder $C:=B_{r_0}^{n-1}(0)\times I_0$, thus yielding a special chart $(U_p,h_p)$ for $p$ with $U_p=h_p^{-1}(C)$.   The resulting atlas $\{(U_p,h_p)\}_{p\in\pa\Om}$ shows that $\Om$ has uniformly Lipschitz boundary.
\end{enumerate}
\end{remark}

\begin{defn}\label{def.unifconv} Let $\{\Om_j\}_{j=1}^\infty$ be a sequence of domains with Lipschitz boundary in a Riemannian manifold $M$, and $\Om\subset M$ be a domain with Lipschitz boundary.  We say the sequence $\{\Om_j\}_{j=1}^\infty$ \emph{converges uniformly} to $\Om$ if
\begin{enumerate}
\item $\lim_{j\to\infty}\, d_H(\pa\Om,\pa\Om_j)=0$ where $d_H$ denotes the Hausdorff distance.
\item \label{def.unifconv_2} A finite collection $\mathcal{C}$ of special charts whose domains cover $\pa\Om$ (see part \ref{part3} of {Remark~\ref{rem.weak}}) can be chosen so they simultaneously cover $\pa\Om_j$ for all $j$ and such that for each $(U,h)\in \mathcal{C}$ with associated data $(r, I, \phi)$ as in Definition~\ref{def.lip}, we have:
\begin{enumerate}
\item $h(U\cap \pa \Om_j)$ is the graph of a Lipschitz function $\phi_j: B_r^{n-1}(0)\to I$ and $$h(U\cap \Om_j)=\{(x,t)\in B_r^{n-1}(0)\times I: t>\phi_j(x)\}.$$
Moreover the Lipschitz constants of the $\phi_j$ are bounded independently of $j$.
\item $\phi_j\to \phi$ uniformly.
\end{enumerate}
\end{enumerate}

\end{defn}

\begin{remark}\label{rem.unif cone} The notion of uniform convergence in Definition~\ref{def.unifconv} differs a priori from the corresponding definition given in \cite{BN20} for sequences of domains with Lipschitz boundary in Euclidean space.  There the second condition is replaced by the uniform cone condition.  But it is known that \eqref{def.unifconv_2} is equivalent to the uniform cone condition; see, e.g., \cite[Theorem 2.4.7]{HP05}.

\end{remark}

The existence of a sequence of smooth domains approximating a domain with Lipschitz boundary is a classical result. The proof is rather technical though, especially when we require this sequence to converge in a \textit{uniform} way. In the Euclidean setting, the result goes back to Jind\v{r}ich Ne\v{c}as \cite{Necas} (see also \cite[Theorem 1.12]{Ver84} and the detailed proof in the appendix to \cite{Ver82}), and it is stated in the manifold setting in \cite[Appendix A]{MT99} where the authors assert that the proof of the Euclidean case can be adapted to the case of Riemannian manifolds.   In dimension two or, more generally, in case all boundary components of $\Om$ are topological spheres, one can give a simple proof in the Riemannian setting by appealing directly to the Euclidean result.   For completeness, we include the general statement (Lemma~\ref{movingdomains}(a,b)) and provide a proof for this special case.

\begin{lemma}\label{movingdomains}  Let $\Om$ be a domain with Lipschitz boundary in a smooth $n$-dimensional Riemannian manifold $M$. Then there exists a sequence of smooth bounded domains $\Om_j$ such that 
\begin{itemize}
    \item[a)] $\Om_j\to \Om$ uniformly in the sense of Definition~\ref{def.unifconv}.
    \item[b)] $\Om_j\subset \Om$.
    \end{itemize}
Moreover, if $\rho\in L^\infty(\pa\Om)$ is a non-negative function, then there exist positive functions $\rho_j\in C^\infty(\pa\Om_j)$ such that:
\begin{itemize}
\item the $\rho_j$ are uniformly bounded independent of $j$;
\item for $\mu_j:=\rho_j \mathcal{H}^{n-1}\vert_{\partial \Om_j}$ and $\mu=\rho {\mathcal{H}^{n-1}} \vert_{\partial \Om}$ viewed as measures on $M$ supported on $\partial\Om_j$ and $\partial\Om$ respectively, we have 
$
\mu_n{\rightharpoonup}  \mu
$
weak-$*$ in the measure sense.  Here $ {\mathcal{H}^{n-1}}$ denotes the {$n-1$}-dimensional Hausdorff measure on $M$.
\end{itemize}
\end{lemma}

\begin{proof}[Proof in the case that the boundary components are topological spheres.] List the connected components of $\pa\Om$ as $\{C_\alpha\}_{\alpha=1}^\ell$.   Since $\Om$ has Lipschitz boundary and $C_\alpha$ is a topological sphere, we may take mutually disjoint open neighborhoods $V_\alpha$ of $C_\alpha$, $\alpha=1,\dots \ell$ and orientation-preserving diffeomorphisms $F_\alpha: V_\alpha\to \R^n$, where each $F_\alpha(V_\alpha)$ is an open set in $\R^n$ diffeomorphic to an annulus.  Moreover, shrinking $V_\alpha$ if necessary, we may assume that the norms $\|dF_\alpha\|$ and $\|dF^{-1}_\alpha\|$ relative to the Riemannian metric on $M$ and the Euclidean metric on $\R^n$ are bounded, so $F_\alpha$ is bi-Lipschitz.  The curve $F_\alpha(C_\alpha)$ bounds a simply connected domain $D_\alpha$ with $D_\alpha\cap F_{\alpha}(V_\alpha)=F_\alpha( \Om\cap V_\alpha)$.   

By Remark \ref{rem.weak}(\ref{part3}), we may choose a collection $\A_\alpha$ of {special} charts whose domains cover $C_\alpha$.  We may assume $U\subset V_\alpha$ for each $(U,h)\in \A_\alpha$.   By {Remark~\ref{rem.weak}(2)}  applied to the charts $(F_\alpha(U), h\circ F_\alpha^{-1})$, we see that $D_\alpha$ is a domain with Lipschitz boundary in $\R^n$.    By \cite[Theorem 1.12]{Ver84}, there exists a sequence of smooth domains $D_\alpha^j\subset D_\alpha$ that converge uniformly to $D_\alpha$ in the sense of Definition~\ref{def.unifconv}.   Moreover, there exist homeomorphisms $\Lambda_j: \pa D_\alpha\to \pa D_\alpha^j$ such that 
\begin{equation}\label{supa}\sup_{x\in \pa D_\alpha}\, |x-\Lambda_j(x)|\,\to 0 \,\mbox{\,as\,\,} j\to \infty.\end{equation}
Finally, there exist positive functions $\omega_\alpha^j\in L^\infty(D_\alpha^j)$ bounded away from zero and infinity uniformly in $j$ such that 
\begin{equation}\label{ome}\Lambda_j^*d\sigma_\alpha^j=\omega_\alpha^j d\sigma_\alpha,\end{equation}
where $d\sigma_\alpha^j$ and $d\sigma$ are the surface measures on $\pa D_\alpha^j$ and $\pa D_\alpha (=F_\alpha( C_\alpha))$, respectively.   

We have $\pa D_{\alpha,j}\subset F{{_\alpha}}(V_\alpha)$ for all sufficiently large $j$ and thus we may assume this is the case for all $j$.  Set
$$C_\alpha^j=F_\alpha^{-1}(\pa D_\alpha^j)$$ 
and let $\Om_j$ be the domain in $M$ bounded by $C_\alpha^1\cup\dots\cup C_\alpha^\ell$.   Then the sequence $\Om_j$ satisfies conditions (a) and (b) of the lemma.   

We now prove the final statement of the lemma. For each $j$, define a homeomorphism 
\begin{equation}\label{ta}\tau_j:\pa\Om \to \pa \Om_j\mbox{\,\, by\,\,} \tau_j|_{C_\alpha} =F_\alpha^{-1}\circ \Lambda_\alpha^j\circ F_\alpha.\end{equation}
By Equation~\eqref{supa}, we have
\begin{equation}\label{supa2}\sup_{p\in \pa \Om}\, |d(p,\tau_j(p))|\,\to 0 \,\mbox{\,as\,\,} j\to \infty\end{equation}
where $d$ denotes Riemannian distance in $M$. Using Equation~\eqref{ome}and the fact that the $F_\alpha$ are bi-Lipschitz, we moreover obtain a sequence of measurable functions $\eta^j:\pa\Om\to \R$ that are uniformly bounded away from zero and infinity such that $\tau_j^*(d\sigma_j)=\eta^j d\sigma$ where $\sigma$ and $\sigma_j$ are the surface measures on $\pa\Om$ and $\pa\Om_j$.  

Now let $\rho\in L^\infty (\pa \Om)$ be non-negative.   Then 
\begin{equation}\label{kaj}(\tau_j^{-1})^*\rho d\sigma=\kappa_j d\sigma_j \mbox{\,\,where\,\,}\kappa_j=
\frac{\rho}{\eta_j}\circ (\tau_j)^{-1}.\end{equation}
The functions $\kappa_j$ are non-negative and $\sup_j \|\kappa_j\|_{L^\infty(\pa\Om_j)} <\infty$.   Since $L(\pa \Om_j)$ is bounded independently of $j$, we also have that $\sup_j \|\kappa_j\|_{L^1(\pa\Om_j)} <\infty$
We now use the fact that $\pa\Om_j$ is a smooth compact manifold to choose a strictly positive function $\rho_j\in C^\infty(\pa\Om_j)$ such that 
\begin{equation}\label{1j}\|\rho_j -\kappa_j\|_ {L^1(\pa\Om_j)} < \frac{1}{j}.\end{equation}
 In particular, we have
 \begin{equation}\label{supa3} \sup_j \|\rho_j\|_{L^1(\pa\Om_j)} <\infty.\end{equation}
 
 To complete the proof, let $f\in C_c^\infty(M)$.    Then by Equation~\eqref{kaj}, we have
 $$\left|\int_{\pa\Om} f\rho d\sigma -\int_{\pa\Om_j} f\rho_j d\sigma_j\right| = \left|\int_{\pa\Om_j}[( f\circ \tau_j^{-1}) \kappa_j- f\rho_j]\,d\sigma_j\right|\leq  A_j +B_j$$ 
 where 
 $$A_j= \left|\int_{\pa\Om_j}\,[( f\circ \tau_j^{-1}) (\kappa_j-\rho_j)\,d\sigma_j\right|\leq (\max|f| )\|\rho_j -\kappa_j\|_ {L^1(\pa\Om_j}$$
 and 
 $$B_j=\left|\int_{\pa\Om_j}\,[( f\circ \tau_j^{-1})-f]\, \rho_j \,d\sigma_j\right|\leq \|( f\circ \tau_j^{-1})-f\|_{L^\infty(\pa\Om_j}\,\|\rho_j\|_{L^1(\pa\Om_j}.$$
 
By Equation~\eqref{1j}, we have $\lim_{j\to\infty}\,A_j=0$.   By uniform continuity of $f$ and Equations~\eqref{supa2} and \eqref{supa3}, we also have $\lim_{j\to\infty}\,A_j=0$, thus completing the proof.\end{proof}

We conclude this background by recalling a recent result of Bucur and Nahon \cite{BN20},  addressing the convergence of Steklov eigenvalues.  Their result was stated in the context of bounded domains in $\R^2$, but it is straightforward to check that the result holds for domains in a  complete Riemannian surface.

\begin{prop}\label{Proposition 2.3, BuNa2.3}
Let $\Om$ be a domain with Lipschitz boundary in a  Riemannian surface $M$ and let $\{\Om_j\}$ be a sequence of  domains in $M$ with Lipschitz boundary.  Also, let $\rho \in L^\infty (\partial \Om)$ be a nonnegative function and $\rho_j \in L^\infty (\partial \Om_j)$ be a sequence of nonnegative functions. Assume that the following conditions are satisfied:
\begin{enumerate}
\item $\Om_j\to \Om$ uniformly in the sense of Definition~\ref{def.unifconv},
\item $
\limsup_{j \to \infty} ||\rho_j||_{L^\infty(\partial \Om_j)} < \infty$\,;
\item for $\mu_j:=\rho_j \mathcal{H}^{1}\vert_{\partial \Om_j}$ and $\mu=\rho \mathcal{H}^{1} \vert_{\partial \Om}$ viewed as measures on $M$ supported on $\partial\Om_j$ and $\partial\Om$ respectively, we have 
$
\mu_j{\rightharpoonup}  \mu
$
weak-$*$ in the measure sense.  
\end{enumerate}   Then for all $k \geq 1$,
\[
\lim_{j \to \infty} \sigma_k(\Om_j, \rho_j) =\sigma_k (\Om, \rho).
\]
\end{prop}

\subsection{Extending eigenvalue bounds from smooth to Lipschitz manifolds}\label{smooth to lipshitz}

\begin{prop} \label{topbound1} Let $(M,g)$ be a complete Riemannian surface.   Suppose that  for all smooth bounded domains $\Om$ in $M$, 
\begin{equation}\label{smoothbd}\sigma_k(\Om,g)L(\partial \Om,g)\leq S,\end{equation}
where $S$ is a positive constant. 
Then inequality~\eqref{smoothbd} holds for all bounded domains $\Om$ in $M$ with Lipschitz boundary.
  Moreover, if $\Om$ is a bounded domain in $M$ with Lipschitz boundary and $\rho\in L^\infty(\partial\Om)$ is any non-negative weight, then we also have
$$\sigma_k(\Om,g,\rho)L_\rho(\Om,g)\leq S.$$

\end{prop}

The proposition is immediate from Lemma~\ref{movingdomains}, Proposition~\ref{Proposition 2.3, BuNa2.3} and Remark~\ref{item:smoothpos}.

\begin{remark}\label{rem.KL}
As we were completing this paper, we learned that Karpukhin and Lagac\'e \cite{KL22} have independently obtained  Proposition~\ref{topbound1} in arbitrary dimension with a  suitable normalisation of eigenvalues.    
\end{remark}

\begin{remark}\label{mixSN}
    As noted in Remark~\ref{rem:linking_the_probs}, a mixed Steklov-Neumann problem can be viewed as a weighted Steklov problem in which the weight $\rho\equiv0$ on the Neumann part of the boundary. Proposition \ref{Proposition 2.3, BuNa2.3} and Lemma \ref{movingdomains} imply that the eigenvalues of mixed Steklov-Neumann problems on domains with Lipschitz boundary can be expressed as limits of Steklov eigenvalues on  {domains with smooth boundary}. 
    
\end{remark}

Consider the class of all Riemannian surfaces of a given topological type.   Within such a class, Proposition~\ref{topbound1} implies that any Steklov eigenvalue bound that is valid for all Riemmanian surfaces with smooth boundary remains valid for all the Riemannian surfaces with Lipschitz boundary.    Moreover, the same bound holds for all weighted Steklov problems with $L^\infty$ non-negative weight on Riemannian surfaces in this class:

\begin{nota}\label{nota:dgb} Let $\gamma\in\{0,1,2\dots\}$ and $b\in \Z^+$. 
Denote by $\mathcal{C}(\gamma,b)$, respectively $\mathcal{L}(\gamma,b)$, the collection of all compact Riemannian surfaces $(\Om,g)$ with $C^\infty$ boundary, respectively Lipschitz boundary, such that
$\Om$ has genus $\gamma$ and $b$ boundary components.   

Define
\begin{eqnarray*}
\Sigma_k(\gamma,b)&:=&\sup\{\sigma_k(\Om, g)L(\partial \Om, g): (\Om_{\gamma,b},g)\in \mathcal{C}(\gamma,b)\},\\
\Sigma_k^*(\gamma,b)&:=&\sup\{\sigma_k(\Om_{\gamma,b}, {\rho}g)L_{\rho}(\partial \Om, g): (\Om,g)\in \mathcal{L}(\gamma,b), 0\le\rho\in L^\infty(\pa\Om)\},
\end{eqnarray*}

\end{nota}

\begin{cor}\label{topbound}
In the notation of~\ref{nota:dgb}, we have
$$\Sigma_k^*(\gamma,b)=\Sigma_k(\gamma,b).$$\
\end{cor}

\subsection{From special to general mixed problems}
The classical sloshing problem is a mixed Steklov-Neumann problem in which the portion of the boundary of the domain with the Steklov boundary condition is a line segment representing the freely moving surface of a fluid; the remainder of the boundary represents the walls of a cross-section of the container, has the Neumann boundary condition, and is a priori unrestricted in terms of its geometry.  In contrast to this setup, we make use of a class of domains with Lipschitz boundary where the Neumann (or Dirichlet) part of the boundary consists of geodesic arcs and the Steklov part of the boundary is less geometrically restricted. We say that these domains satisfy the \emph{Lipschitz doubling condition}, and fully define them below.  {The subsequent Lemmas~\ref{lem.doub lip} and \ref{remove corners} together with Lemma \ref{approx by polygonal} and Proposition~\ref{doub to general} in the next section will enable us} to prove eigenvalue bounds and asymptotics for much more general classes of mixed Steklov-Neumann and Steklov-Dirichlet problems on surfaces in Sections \ref{bounds_mixed} and \ref{sec:4}.   A second motivation comes from interest in orbifolds as explained in Remark~\ref{nota.orb} below.   We will not use the language of orbifolds in this article, and knowledge of orbifolds is not needed for any of our results.

\begin{defn}\label{def.doub} We will say that  $(\Om,g)$   with Lipschitz boundary, together with a non-trivial decomposition $\partial\Om=\psom \sqcup \partial_*\Om$, satisfies the \emph{Lipschitz doubling condition} if the following conditions hold:
\begin{enumerate}
\item Each component of $\partial_*\Om$ is either a non-trivial geodesic segment (meeting $\psom$ only at its endpoints) or a closed geodesic; 
\item At each point $p$ where $\psom$ and $\partial_*\Om$ meet, there exists a neighborhood $U$ of $p$ in the ambient Riemannian surface $M$ such that $U\cap \Om$ is contained in a geodesic wedge emanating from $p$ of angle strictly less that $\pi$.    (If the boundary is sufficiently regular so that the interior angle $\alpha$ between $\psom$ and $\partial_*\Om$ is well-defined, this condition just says that $0<\alpha<\pi$.)

\end{enumerate}
Neumann, Dirichlet, or mixed Dirichlet-Neumann boundary conditions will be placed on $\partial_*\Om$.
\end{defn}

\begin{lemma}\label{lem.doub lip} Suppose that $\Om$, together with the decomposition $\partial\Om=\psom \sqcup \partial_*\Om$, satisfies the Lipschitz doubling condition.   Then by doubling $\Om$ across $\partial_*\Om$ we obtain a Riemannian surface $\Sigma$ with Lipschitz boundary.

\end{lemma}

\begin{proof}
We first show that the double embeds in a complete Riemannian manifold.
The fact that the double is a Riemannian manifold with Lipschitz boundary is immediate from the doubling condition.  We show that the double $\Sigma$ satisfies Definition \ref{def.lip}.\ref{stand}.   

By the standing assumption, $\Om$ is a domain in some complete Riemannian manifold $M$.   We will perturb $\Om$ to obtain a bounded domain $\Om'$ in $M$ and a boundary decomposition $\partial\Om'=\psom \cup \partial_*\Om'$ in such a way that:
\begin{enumerate}
\item $\Om \subset \Om'$
    \item $\partial_* \Om'\cap \Om =\partial_* \Om$.
 \item Each component of $\psom'$ is smooth.
 \item At any point where $\partial_* \Om'$ and $\psom'$ meet, they are orthogonal.
 \item $\Om'$ satisfies the Lipschitz doubling condition.
\end{enumerate}

The double $\Sigma'$ of $\Om'$ across $\partial_* \Om'$ will then be a compact Riemannian manifold containing $\Sigma$.  Moreover, the last three conditions imply that $\Sigma'$ has smooth boundary and thus embeds smoothly as a domain in a complete Riemannian manifold $M'$.  Since $\Sigma\subset\Sigma'$, the standing assumption holds for $\Sigma$.  

To obtain a domain $\Om'$ in $M$ satisfying these conditions, perturb $\partial\Om$ as follows:  
\begin{itemize}
    \item Any component of $\partial \Om$ that lies wholly in $\partial_*\Om$ is left unchanged and will form part of $\partial_*\Om'$.
    \item Any component of $\partial \Om$ that lies wholly in $\psom$ is continuously perturbed in the exterior of $\Om$ to obtain a smooth simple closed curve that will form a component of $\psom'$.
    \item We perturb each component  $C$ of $\partial \Om$ that intersects both $\partial_*\Om$ and $\psom$ in the following way:  Each component $\gamma$ of $C\cap \partial_*\Om$ is a segment of a geodesic in $M$.  The Lipschitz doubling condition guarantees that this geodesic exits $\Om$ when one reaches either endpoint of $\gamma$.  Extend each such geodesic segment slightly from both endpoints.   We can then perturb the segments of $C\cap \psom$ into the exterior of $\Om$ to form smooth segments that meet the extended geodesic segments at right angles.
\end{itemize}

Carry out this process one component at a time, being careful at each step that the newly perturbed boundary components stay away from the previously perturbed ones.
\end{proof}

\begin{lemma}\label{remove corners} Let $(\Om,g)$ be a compact Riemannian surface with Lipschitz boundary $\partial\Om= \psom \sqcup\partial_*\Om$ and suppose that each connected component of $\partial_*\Om$ is polygonal (i.e., piecewise geodesic).    Let $\{q_1,\dots, q_n\}$ be the set of all the corners of all the polygonal paths making up $\partial_*\Om$.   (We do not include in this list any endpoints where $\partial_*\Om$ meets $\psom$.)  Then there exists a Riemannian metric $g'$ on $\Om$ such that:
\begin{enumerate}
\item $g'$ coincides with $g$ on a neighborhood of $\psom$;
\item $g'$ is conformally equivalent to $g$ on $\Om\setminus\{q_1,\dots q_n\}$.   (At each $q$, the conformal factor approaches zero or infinity depending on the angle.)
\item Each component of $\pst\Om$ is a geodesic segment or a closed geodesic.
\end{enumerate}
Consequently, the metrics $g$ and $g'$ on $\Om$ have both the same Steklov-Neumann spectrum and the same Steklov-Dirichlet spectrum, where we place Steklov boundary conditions on $\psom$ and Neumann or Dirichlet conditions on $\partial_*\Om$.   
\end{lemma}

Note that if $(\Om,g)$ satisfies the second condition in Definition~\ref{def.doub}, then $(\Om,g')$ satisfies the Lipschitz doubling condition.

\begin{remark} To clarify the statement of the lemma and the proof below, the definition~\ref{stand} of Riemannian manifold with Lipschitz boundary says that $(\Om,g)$ is a bounded domain in some Riemannian manifold $M$.   However, if we remove a neighborhood of $\psom$, then -- ignoring the ambient manifold  -- the remainder of $(\Om,g)$ can be viewed in various additional ways: (i) as a manifold with corners along with a smooth Riemannian metric or (ii) as a  manifold with smooth boundary along with a Riemannian metric that has conical singularities at the points $q_j$ on the boundary.  Taking the latter point of view, we will construct the metric $g'$ by making a conformal change in $g$ that removes the conical singularities.  With this new metric, we will see that $(\Om, g')$ will again satisfy the definition of a Riemannian manifold with Lipschitz boundary; however it can be realized as a domain in a \emph{different} Riemannian manifold $(M',g')$.  The metric $g'$ does \emph{not} extend smoothly to a Riemannian metric on the original ambient manifold $M$.  Any such extension will have singularities at the points $q_j$.

\end{remark}

\begin{proof}
Fix $j\in \{1,\dots, n\}$.  Let $c\pi\in (0,2\pi)$ be the interior angle at $q$.    In geodesic polar coordinates on a sufficiently small neighborhood $U$ of $q$ in $\Om$, the Riemannian metric has the form
$$g=ds^2+\eta(s,\tau)d\tau^2, \hspace{1cm}(s,\tau)\in [0,\epsilon)\times [0,c\pi].$$ 
The fact that $g$ is of class at least $C^2$ on the ambient space says that 
  \begin{equation}\label{smmet}\lim_{s\to 0}\frac{1}{s^2} \eta(s,\tau) = 1 , \hspace{1cm} \lim_{s\to 0}\frac{1}{2s}\frac{\partial}{\partial s} \eta(s,\tau) = 1 , \hspace{1cm} 
\lim_{s\to 0}\frac{1}{2}\frac{\partial^2}{\partial s^2} \eta(s,\tau)= 1.\end{equation}

Making the change of variables $\theta=\frac{1}{c}\tau$, we have 
$$g=ds^2+c^2\kappa(s,\theta)d\theta^2, \hspace{1cm}0\leq s<\epsilon, \,\,\,0\leq \theta\leq \pi$$
on $U$, 
where $\kappa(s,\theta)=h(s,c\theta)$.    Observe that $\kappa$ satisfies the same conditions~\eqref{smmet} as $\eta$.

Next making the change of variables $s=r^{c}$, we obtain
$$g=\alpha(r)\left(dr^2+\lambda(r,\theta)d\theta^2\right)=\alpha(r)h$$
with
\begin{equation}\label{eqn.alpha}
\alpha(r)=c^2 r^{2(c-1)},\hspace{1cm}\lambda(r,\theta)=r^{2(1-c)}\kappa\left(r^{c},\theta\right),\,\hspace{1cm} h=dr^2+\lambda(r,\theta)d\theta^2 .
\end{equation}
Straightforward computation confirms that as $r \rightarrow 0$ we have 
$$\frac{1}{r^2} \lambda(r,\theta) \rightarrow 1, \hspace{1cm} \frac{1}{2r}\frac{\partial}{\partial r} \lambda((r,\theta) \rightarrow 1,\hspace{1cm}\frac{1}{2}\frac{\partial^2}{\partial r^2} \lambda((r,\theta) \rightarrow 1.$$ 
Thus $h$ is a $C^2$-smooth metric and is conformal to $g$ on $U\setminus\{q\}$.
The conformal factor has limit either 0 or $\infty$ at $q$, depending on whether the original interior angle $c\pi$ is less than or greater than $\pi$.   Moreover, relative to the metric $h$, the intersection $U\cap \partial\Om$ is a geodesic arc.  

Note that $U$ is given in our coordinates $(r,\theta)$ by $\{r<\delta\}$ (where $\delta=\epsilon^{1/c}$).   Let $0<\delta_1<\delta_2<\delta$ and define subsets $V_1$ and $V_2$ of $U$ by
$$
V_1 = \{ (r, \theta): r < \delta_1 \}
\ \ \ \text{and} \ \ \ 
V_2 = \{ (r, \theta): r < \delta_2 \}.
$$
Let $\rho \in C^\infty(0,\delta)$ satisfy
$$\rho(r)=\begin{cases}1&\hspace{1cm}0\leq r\leq\delta_1\\\alpha(r)& \hspace{1cm}\delta_2\leq r<\delta\end{cases}$$
where $\alpha(r)$ is as in \eqref{eqn.alpha}.
Then the smooth Riemannian metric $h'=\rho(r)h$ coincides with $h$ on $V_1$ and with $g$ on $U\setminus{V_2}$.  Moreover, since the conformal factor depends only on $r$, the intersection of $U$ with $\partial\Om$ is still geodesic.   

Carry out the procedure above for each $j$ to obtain a metric $h'_j$ as above on a neighborhood $U_j$ of $q_j$, with the $U_j$ chosen sufficiently small that they are mutually disjoint and do not intersect $\psom$.    Set $g'=h'_j$ on $U_j$, $j=1,\dots, n$ and $g'=g$ on $\Om\setminus(U_1\cup\dots\cup U_n)$.  Then ${g'}$ satisfies the first two itemized conditions in the lemma.

It is straightforward to construct a complete Riemannian manifold $(M',g')$ so that $(\Om,g')$ can be realized as a domain with Lipschitz boundary in $(M', g')$.  

 The final statement of the lemma follows from the variational characterization of the eigenvalues; the fact that the conformal factor is singular on a set of measure zero does not affect the Rayleigh quotient.
 \end{proof}
{\begin{remark}\label{add corners}
Very minor adjustments to the argument in Lemma~\ref{remove corners} show that if some open segment $S$ of the boundary of a  compact Riemannian surface $(\om,g)$ is polygonal, then one can conformally adjust the interior angle at any point $q$ of $S$ by making a conformal change $g'$ of metric with $g'=g$ outside an arbitrarily small neighborhood of $q$ in $\om$.  In particular, this allows us to ``add corners'' into a geodesic boundary segment via a conformal change of metric.

\end{remark}}

We end this section with an elementary but important example.

\begin{ex}[Model examples]\label{models}\begin{itemize}
    \item [(i)] Let $\disk(\ell)$ denote the Euclidean disk of circumference $\ell$. Then 
    \[
\operatorname{Stek}(\disk(\ell)) = \left\{0, \frac{2\pi}{\ell}, \frac{2\pi}{\ell}, \frac{4\pi}{\ell}, \frac{4\pi}{\ell}, \frac{6\pi}{\ell}, \frac{6\pi}{\ell}, {\dots}\right\},
\]
viewing the right-hand-side as a multiset. Moreover, the Steklov eigenfunctions are $r^k \cos(kt)$ and $r^k \sin (kt)$. 
\item[(ii)] Let $\tau$ be the reflection across a diameter of $\disk(\ell)$. The quotient of the Euclidan disk $\disk(\ell)$ by the reflection $\tau$ can be identified with a flat half disk. We denote it by $\hD(\bar\ell)$, where $2\bar\ell=\ell$. Let consider the boundary decomposition $\pa\hD(\bar\ell)=\pa_S\hD(\bar\ell)\sqcup\pa_*\hD(\bar\ell)$ where $\pa_S\hD(\bell)$ is a  half-circle and $\pa_*\hD(\bar\ell)$ is the diameter across which we reflected, i.e. the fix point of $\tau$. Then 
\[
\Stek_N(\hD(\bell)) = \left\{0, \frac{\pi}{\bell}, \frac{2\pi}{\bell},  \frac{3\pi}{\bell}, {\dots}\right\},\qquad
\Stek_D(\hD(\bell)) = \left\{ \frac{\pi}{\bell}, \frac{2\pi}{\bell},  \frac{3\pi}{\bell}, {\dots}\right\}.
\]
Notice that $\sigma_k^N(\hD(\bell))=\sigma_{k-1}^D(\hD(\bell))=\frac{\pi k}{\bell}$, $k\ge1$.
\item[(iii)] Let $\qD(\bell)$ denote the quotient of $\disk(4\bell)$ by two orthogonal reflections $\tau_1$ and $\tau_2$. We identify $\qD(\bell)$ with a quarter of a flat disk. Let consider the boundary decomposition $\pa\qD(\bar\ell)=\pa_S\qD(\bar\ell)\sqcup\pa_*\qD(\bar\ell)$ where $\pa_S\qD(\bell)$ is a  quarter-circle.  We also have
\begin{equation}\label{qD vs hD}\Stek_D(\qD(\bell))=\Stek_D(\hD(\bell)),\quad \Stek_N(\qD(\bell))=\Stek_N(\hD(\bell)). \end{equation}
Note that the radius of $\qD(\bell)$ is twice the radius of $\hD(\bell)$.\\
Let us consider $\pa_*\qD(\bar\ell)=\pa_N\qD(\bell)\sqcup\pa_D\qD(\bell)$, where $\pa_N\qD(\bell)$ is the straight segment fixed by $\tau_1$ and $\pa_D\qD(\bell)$ is the straight segment fixed by $\tau_2$.   We impose the Neumann boundary condition on $\pa_N\qD(\bell)$ and  Dirichlet boundary condition on $\pa_D\qD(\bell)$.  The Spectrum of this mixed Steklov problem $\Stek_{\rm DN}(\qD(\bell))$ is given by 
\[\Stek_{\mix}(\qD(\bell))=\frac{\pi}{2\bell}\left\{1,3,5,\cdots\right\}.\]
Indeed, if we view  $\hD(2\bell)$ obtained by the quotient of $\disk(4\bell)$ by reflection $\tau_1$, then  $\qD(\bell)$ can be viewed as the quotient of $\hD(2\bell)$ by a reflection $\tau_2$. The Steklov-Dirichlet eigenfunctions on $\Stek_D(\hD(2\bell))$ can be grouped into invariant and anti-invariant with respect to $\tau_2$. Hence,
\begin{equation}\label{hD vs qD+qD}\Stek_D(\hD(2\bell))=\Stek_{\mix}(\qD(\bell))\sqcup \Stek_D(\qD(\bell)).\end{equation}
Hence, $\Stek_{\mix}(\qD(\bell))=\Stek_D(\hD(2\bell))\setminus\Stek_D(\qD(\bell))=\frac{\pi}{2\bell}\left\{1,3,5,\cdots\right\}.$

\end{itemize}
We denote the unit disk by $\disk$ instead of $\disk(2\pi)$ for simplicity. Similarly in parts (ii) and (iii), we denote $\hD(\pi)$ and $\qD(\pi/2)$  by $\hD$ and $\qD$ since they can be viewed as quotient of the unit disk.
\end{ex}

For a general Riemannian surface with boundary with reflection symmetry we introduce the following notation that we use in the following sections. 

\begin{notarem}\label{nota.orb}
\begin{itemize}\item[(i)] Let $\Om$ be a Riemannian surface with boundary and suppose that $\Om$ admits a reflection symmetry $\tau$.   We will denote the quotient surface $\Om/\langle\tau\rangle$ by $\Om_{\tau}$.    The quotient can itself be viewed as a Riemannian surface, with $\pa\Om_\tau =\pa_1\Om_\tau \cup \pa_2\Om_\tau $, where $\pa_1\Om_\tau$ is the image of $\pa\Om$ and $\pa_2\Om$ is the image of the fixed point set of $\tau$.  The Steklov-Neumann, respectively Steklov-Dirichlet, problem on $\Om_\tau$ with Neumann, respectively Dirichlet, with $\pa_S\Om_\tau=\pa_1\Om_\tau$ and Neumann, respectively Dirichlet, conditions on $\pa_*\Om_\tau=\pa_2\Om_\tau$  coincides with the part of the Steklov spectrum of $\Om$ for which the eigenfunctions are $\tau$-invariant, respectively $\tau$-anti-invariant.  Surfaces satisfying the Lipschitz doubling condition are precisely the surfaces that can be realized as quotients $\Om_\tau$ for which both $\pa\Om$ and the fixed point set of $\tau$ are connected.
\item[(ii)]  Every Riemannian orbisurface with non-trivial boundary can be identified, up to a (possibly singular) conformal change of metric, with a quotient surface $\Om_\tau$ as in (i).  When viewed as an orbifold, $\pa_2\Om_\tau$ is viewed not as part of the orbifold boundary but rather as the singular set of the orbifold.   The Steklov spectrum of the orbifold is precisely the mixed Steklov-Neumann spectrum of $\Om_\tau$ with boundary conditions as in (i).  We refer to \cite{ADGHRS,gordon12} for details on orbisurfaces and the Steklov problem on them.  

\end{itemize}

\end{notarem}

\section{Eigenvalue upper bounds for mixed Steklov problems}
\label{bounds_mixed}

\begin{nota}\label{nota.mixed} Let $(\Om,g)$ be a compact orientable Riemannian surface with Lipschitz boundary and a nontrivial decomposition  $\pa\Om=\psom\sqcup\pa_*\Om$.   We denote by  $\sigma_k^D(\Om,\psom)$, respectively $\sigma_k^N(\Om,\psom)$, the $k$th non-zero eigenvalue of the Steklov-Dirichlet, respectively Steklov-Neumann, problem on $\Om$ where Steklov conditions are placed on $\psom$ and Dirichlet, respectively Neumann, conditions are placed on $\pa_*\Om$.    Normally we will suppress the dependence on the decomposition $\pa\Om=\psom\sqcup\pa_*\Om$ and write simply $\sigma_k^D(\Om)$, respectively $\sigma_k^N(\Om)$.

We will say that the boundary decomposition is of type $(b; r; m_1,\dots m_{b-r})$ if $\pa\Om$ consists of $b$ components, exactly $r$ of which lie entirely in $\psom$; each of the remaining $b-r$ boundary components  contains $m_i$ disjoint arcs of $\psom$, $i=1, \dots, b-r$, where $m_i$ may equal $0$.    Here $r\in \{0,\dots,b-1\}$ and $m_1,\dots, m_{b-r}\in\{0,1,2,\dots\}$.   (The only constraint is that $r+m_1+\dots m_{b-r}>0$ since $\psom\neq\emptyset$.)

\end{nota}
\begin{thm}\label{general minpi} In the notation of~\ref{nota.mixed} and \ref{nota:dgb}, suppose that $\Om$ has genus zero and that $\pa\Om$ has boundary decomposition of type $(b;b-1;m)$ for some $m$.  Then we have
\begin{equation}\label{gen min}\min\{\sigma_{k}^N(\Om), \sigma_{k-1}^D(\Om)\}L(\psom)\leq \frac{1}{2}\Sigma_{2k-1}(0,2(b-1)+m),\qquad  k\ge1.
\end{equation}

Moreover, if equality holds, then $\sigma_{k-1}^D(\Om)=\sigma_k^N(\Om)$. 
\end{thm}
It is known (see \cite{KS,GKL,Kok14}) that 
\begin{equation}\label{inq:gkl}
\Sigma_k(0,b)<8\pi k.    
\end{equation} Inequality\eqref{inq:gkl} is sharp but the number of boundary components of any sequence of domains approaching $8\pi k$ cannot remain bounded.  As a result, we have the following immediate corollary. 

\begin{cor}\label{minpiwithholes}
Under the assumption of Theorem \ref{general minpi} we have

\begin{equation}\label{min2}\min\{\sigma_k^N(\Om), \sigma_{k-1}^D(\Om)\}L(\psom)< 4(2k-1)\pi, \qquad  k\ge1.\end{equation}

\end{cor}

We will first prove Theorem~\ref{general minpi} and note a corollary of the proof of possible independent interest. 
We will then address several special cases of the theorem in which stronger conclusions hold.

\subsection{Proof of Theorem~\ref{general minpi}.}
\begin{lemma}\label{DoMono}  Let $\Om_1,\Om_2\in\Lc(\gamma,b)$ (see Notation~\ref{nota:dgb}) with $\partial\Om_i= \psom_i \sqcup\partial_*\Om_i$, $i=1,2$. Suppose that $\Om_1$ is a proper subdomain of $\Om_2$ satisfying $\psom_1=\psom_2$. Then we have
$\sigma_k^D(\Om_1) > \sigma_{k}^D(\Om_2)$ and
$\sigma_k^N(\Om_1) < \sigma_{k}^N(\Om_2)$ for all $k\geq 1$.  
\end{lemma}
This lemma is a generalization of Propositions 3.1.1 and 3.2.1 in \cite{BKPS} and the same proof is valid.

\begin{proof}[Proof of equality statement in Theorem~\ref{general minpi}]
Let  $A:=\frac{1}{2}\Sigma_{2k-1}(0,2(b-1)+m)$. First suppose that $\sigma_{k-1}^D(\Om)L(\psom)=A$ and $\sigma_k^N(\Om)L(\psom)>A$.   Perturb $\pa_*\Om$, leaving $\psom$ unchanged, to obtain a new domain $\Om_0$ properly contained in $\Om$.  By Lemma~\ref{DoMono}, we have  $\sigma_{k-1}^D(\Om_0)>\sigma_{k-1}^D(\Om)$ and $\sigma_k^N(\Om_0)<\sigma_k^N(\Om)$.   In particular, a small enough perturbation of this type will lead to a domain $\Om_0$ for which inequality~(\ref{gen min}) fails.  Similarly, if $\sigma_k^N(\Om)L(\psom)=A$ and $\sigma_{k-1}^D(\Om)L(\psom)>A$, then a suitably chosen small perturbation of $\pa_*\Om$ in the ambient manifold leads to a domain properly containing $\Om$ for which  inequality \eqref{gen min} fails.   Thus we obtain a contradiction in either case.

\end{proof}

\begin{lemma}\label{lem.minpi doub case} Theorem~\ref{general minpi} holds for all $(\Om,g)$ with $\psom\sqcup\pa_*\Om$ that satisfy the Lipschitz doubling condition.  (See Definition~\ref{def.doub}.)

\end{lemma}
\begin{proof}
Double $\Om$ across $\pa_*\Om$ to obtain a smooth surface $\Sigma$ with Lipschitz boundary.    The surface $\Sigma$ has genus zero and $2(b-1)+m$ boundary components.  $ \Sigma$ admits a reflection symmetry $\tau$, and the quotient of $\Sigma$ by $\tau$ is isometric to $\Om$. Every Steklov eigenspace of $\Sigma$ decomposes into $\tau$-invariant and $\tau$-anti-invariant subspaces. The former project to the Steklov-Neumann eigenfunctions and the latter to the Steklov-Dirichlet eigenfunctions on $\Om$. Hence, one of $\sigma_k^N(\Om)$ or $\sigma_{k-1}^D(\Om)$ must appear as $\sigma_q(\Sigma)$ for some $q\leq 2k-1$. Inequality~(\ref{gen min}) immediately follows.

\end{proof}

\begin{lem}\label{approx by polygonal} Let $(\Om,g)\in\Lc(\gamma,b)$ with $\partial\Om= \psom \sqcup\partial_*\Om$.   Then there exists a sequence of bounded domains $\{\Om_j\}$ in $\Om$ such that: 
\begin{enumerate}
\item  $\pa\Om_j=\psom_j\sqcup\pa_*\Om_j$ with $\pa_S\Om_j=\pa_S\Om$.
\item $\pa\Om_j$ satisfies the second condition in Definition~\ref{def.doub};
\item  Each component of $\partial_*\Om_j$ is polygonal;
\item $\Om_j\to \Om$ uniformly in the sense of Definition~\ref{def.unifconv}.  (In particular, $\pa_*\Om_j \to \pa_*\Om$.)
\end{enumerate}

\end{lem}

\begin{proof} It is straightforward to construct a sequence satisfying the first, third and fourth conditions in the lemma.  Thus we will start with such a sequence $\{\Om_j\}$  and then perturb $\partial_*\Om_j$ in a small neighborhood of any endpoint of $\psom$ where the second condition in Definition~\ref{def.doub} fails to hold in order to get a new sequence that satisfies all four conditions.  If there is more than one point where the condition fails, we take the neighborhood around each one small enough so that they are mutually disjoint.
Let $p$ be such a point where the second condition fails. \smallskip  

\noindent\underline{Special case.} For simplicity, we first consider the case that $\Om$ is a domain in $\R^2$.  After rotating and translating, we may also assume that $p$ coincides with the origin $0$ and, within a rectangular neighborhood $U$ of the origin, $\pa\Om$ is the graph of a Lipschitz function $y=\varphi(x)$.  
\begin{figure}
    \centering
  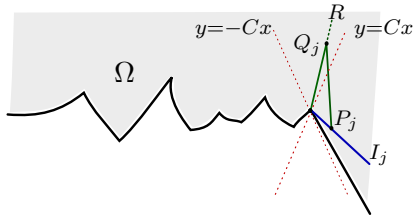
    \caption{An illustration of perturbing $\Om_j$ in Special case of the proof of Lemma \ref{approx by polygonal}.}
   \label{fig:domain0}
\end{figure}
Figure~\ref{fig:domain0} depicts the curve $y=\varphi(x)$; the shaded area above this boundary curve lies in $\Om$.  The part of the curve to the left, respectively right, of the $y$-axis lies in $\psom$, respectively $\pa_*\Om$.   The lines $y=\pm Cx$, where $C>0$ is a Lipschitz constant for $\varphi$ at $0$, appear as dashed lines in Figure~\ref{fig:domain0}.    We may assume for all $j$ that the polygonal curve $(\pa_*\Om_j)\cap U$ lies in the region $\varphi(x)\leq y <Cx$ in the first and fourth quadrants.   We denote by $I_j$ the initial segment, starting from the origin, of the polygonal curve ${\pa_*\Om_j}$.   A sample choice of $I_j$ is depicted in Figure~\ref{fig:domain0}.

We fix a ray $R$ emanating from the origin in the first quadrant lying between the $y$-axis and the line $y=Cx$ as in Figure~\ref{fig:domain0}.    To achieve the third condition, we will perturb the $\Om_j$'s so that their intersections with $U$ lie in the wedge starting at $R$ and continuing counterclockwise to the ray $y=Cx$, $x\leq 0$.  Observe that this wedge has angle less that $\pi$.

Fix a constant $A>C$.    For each $j$, choose a point $P_j\in I_j$ with $|P_j|<\frac{1}{j}$ sufficiently small so that the line through $P_j$ of slope $-A$ intersects $R$ at a point $Q_j\in  R\cap U$.  Now perturb $\Om_j$ by replacing the interval $0P_j$ of $I_j$ by the line segment $0Q_j$ followed by $P_jQ_j$ as in Figure~\ref{fig:domain0}.  Note that the new angles introduced in this polygonal path are uniformly bounded independently of $j$.

After carrying out the procedure above as needed at the various endpoints of $\psom$, one easily checks that the sequence $\{\Om_j\}$ of perturbed domains satisfies all four conditions of the lemma.\smallskip 

\noindent\underline{General case.}  $\Om$ is a bounded domain with Lipschitz boundary.   By Remark~\ref{rem.weak}, we can choose a coordinate chart $(U,\log)$ on a neighborhood $V$ of $p$ in $M$ that is simultaneously a convex geodesic normal chart and a special chart for $\Om$ in the sense of Definition~\ref{def.lip}.   Write $U=\log(V)\subset T_p(M)\simeq \R^2$.  We can then perturb the domains $\log(\Om_j\cap U)\subset \R^2$ as in the special case above (with one modification explained below) and then go back to $U$ via $\exp_p$.    

We clarify the needed modification.   We use the same notation as in the special case; e.g., $I_j$ is the image under $\log$ of the initial geodesic segment issuing from $p$ of $\pa_*\Om$ and $Q_j$ is chosen as above on $I_j$, etc.    Observe that $\exp_p$ carries the straight line segment $0P_j$ to a geodesic segment issuing from $p$ in $\Om_j$.  However, the image under $\exp_p$ of the line segment $P_jQ_j$ will not be a geodesic in general.  Thus we replace this segment by the unique curve $\gamma_j$ from $P_j$ to $Q_j$ such that $\exp_p(\gamma_j)$ is a geodesic.    Note that as $j\to\infty$, we have $P_j\to0$ and $Q_j\to 0$.  By the behavior of $\exp_p$ near the origin, after taking a tail of the sequence, we can thus uniformly bound the angles between $\gamma_j$ and the adjacent segments of the polygonal path.   It is now straightforward the new perturbed sequence satisfied all desired conditions and it    completes the proof.

\end{proof}

\begin{proof}[Completion of Proof of Theorem~\ref{general minpi}]
Given $(\Om,g)$ with  $\pa\Om=\psom\sqcup\pa_*\Om$ as in the theorem, let $\{\Om_j\}\subset \Om$ be a sequence of bounded domains satisfying Lemma \ref{approx by polygonal}.  By Lemmas~\ref{remove corners} and \ref{lem.minpi doub case}, Inequality~\eqref{gen min} holds for each $\Om_j$.   We can view the mixed Steklov-Neumann problem on $\{\Om_j\}$ as a weighted Steklov problem in which the weight $\rho_j$ takes on only the values one and zero on $\psom_j$ and $\pa_*\Om_j$, respectively. Then by Proposition \ref{Proposition 2.3, BuNa2.3}, we get 
\begin{equation}\label{buna}\lim_{j\to\infty}\sigma_k^N(\Om_j) =\sigma_k^N(\Om),\end{equation}
and by Lemma \ref{DoMono} we have $\sigma^D_{k-1}(\Om)\leq\sigma^D_{k-1}(\Om_j)$. Thus Inequality~\eqref{gen min} holds for $(\Om,g,\psom,\pa_*\Om)$ and the theorem follows.

\end{proof}

  The following is a consequence of the proof of Theorem~\ref{general minpi}:
\begin{cor}\label{doub to general}
Fix a boundary decomposition type $T=(b; r; m_1,\dots m_{b-r})$ and fix $\gamma\in \{0,1,2,\dots\}$.   For $\Lc(\gamma,b)$ as in Notation~\ref{nota:dgb}, let $\Lc(\gamma,b, T)$ denote the collection of all $(\Om,g)\in \Lc(\gamma,b)$ along with boundary decompositions $\pa\Om=\psom\cup\pa_*\Om$ of type $T$, and let $\Lc_{\operatorname{doub}}(\gamma,b, T)\subset \Lc(\gamma,b, T)$ denote the subcollection of those that are Lipschitz doublable.    Then any eigenvalue bound for $\sigma_k^N(\Om,g)L(\psom)$ that is valid for all elements of $\Lc_{\operatorname{doub}}(\gamma,b, T)$ is also valid for all elements of $\Lc(\gamma,b, T)$.  The same conclusion holds for mixed Steklov-Dirichlet eigenvalue bounds. 

\end{cor}

\subsection{Special cases of Theorem~\ref{general minpi}.}

We continue to assume that $\pa\Om$ has a boundary decomposition of type $(b;b-1;m)$ as in Theorem \ref{general minpi}. 
 We now address three pairs $(b,m)$  for which there exists a Riemannian metric $g$ achieving equality in Theorem\ref{general minpi} when $k=1$, i.e.,  
\begin{equation}\label{minpiequal}
\sigma_1^N(\Om,g)L(\psom)= \sigma_0^D(\Om,g)\}L(\psom)= \frac{1}{2}\Sigma_{1}(0,2(b-1)+m).
\end{equation}
   First recall that Fraser and Schoen \cite{FS2} showed that $\Sigma_1(0,2)$ is realized by the free boundary minimal surface $(\Om^*,g^*)$ in the Euclidean unit ball $\B^3$  given by the intersection of $\B^3$ with the so-called critical catenoid. (For convenience, we will refer to $(\Om^*,g^*)$ itself as the critical catenoid.)   Moreover, they used this fact to obtain the approximation $\Sigma_1^*(0,2)\sim\frac{4\pi}{1.2}$.  The critical catenoid is centered at the origin and is invariant under reflection  across the equatorial plane and also under reflection  across a plane through its axis. The argument in the proof of Lemma~\ref{lem.minpi doub case} immediately yields:

\begin{prop}\label{special cases} Let $k=1$. Equality is achieved in Theorem \ref{general minpi} in the following cases: 
\begin{enumerate}
    \item $(b,m)=(1,1)$.   The flat half disk in Example~\ref{models} realizes equality.
    \item $(b,m)=(2,0)$  Equality is achieved by the intersection of the critical catenoid $(\Om^*,g^*)$ with the closed upper half plane.
    \item $(b,m)=(1,2)$.  Equality is achieved by the half of $(\Om^*,g^*)$ lying on one side of a plane through its axis, together with its boundary.
\end{enumerate}
In the latter two cases, $\partial_*\Om$ consists of the part of the boundary of the half-catenoid that lies in the cutting plane.
\end{prop}

\begin{thm}\label{minpi}\label{minpisharp} Let $\Omega$ be a simply-connected compact surface with Lipschitz boundary and $\partial \Om= \psom \sqcup \pa_*\Om$. Assume that $\pa\Om$ is of type $(1;0;1)$ as in Notation \ref{nota.mixed}.  Then we have the sharp inequality:
 \begin{equation}\label{min}\min\{\sigma_k^N(\Om), \sigma_{k-1}^D(\Om)\}L(\psom)\leq (2k-1)\pi.\end{equation}

\end{thm}

\begin{proof}
The inequality follows from Theorem~\ref{general minpi} and the Hersch-Payne-Schiffer bound $\Sigma_j(0,b)\leq 2\pi j$.    Girouard and Polterovich \cite{GP} proved that the Hersch-Payne-Schiffer bound is sharp by constructing a collection of domains $\tilde\Om_\epsilon$, each consisting of $2j-1$ overlapping disks as in Figure \ref{gpexample} such that 
$$\lim_{\epsilon\to 0}\,\sigma_j(\tilde\Om_\epsilon)L(\pa\tilde\Om_\epsilon)=2\pi j.$$   
Each of these domains admits a reflection symmetry about the vertical axis of the central disk.  Letting $\Om_\epsilon$ denote the intersection of $\tilde\Om_\epsilon$ with the closed half plane bounded by this vertical axis, we can follow the argument in Lemma~\ref{lem.minpi doub case} to complete the proof.
\begin{figure}
    \centering
\begingroup%
  \makeatletter%
  \providecommand\color[2][]{%
    \errmessage{(Inkscape) Color is used for the text in Inkscape, but the package 'color.sty' is not loaded}%
    \renewcommand\color[2][]{}%
  }%
  \providecommand\transparent[1]{%
    \errmessage{(Inkscape) Transparency is used (non-zero) for the text in Inkscape, but the package 'transparent.sty' is not loaded}%
    \renewcommand\transparent[1]{}%
  }%
  \providecommand\rotatebox[2]{#2}%
  \newcommand*\fsize{\dimexpr\f@size pt\relax}%
  \newcommand*\lineheight[1]{\fontsize{\fsize}{#1\fsize}\selectfont}%
  \ifx\svgwidth\undefined%
    \setlength{\unitlength}{278.23685245bp}%
    \ifx\svgscale\undefined%
      \relax%
    \else%
      \setlength{\unitlength}{\unitlength * \real{\svgscale}}%
    \fi%
  \else%
    \setlength{\unitlength}{\svgwidth}%
  \fi%
  \global\let\svgwidth\undefined%
  \global\let\svgscale\undefined%
  \makeatother%
  \begin{picture}(1,0.30593835)%
    \lineheight{1}%
    \setlength\tabcolsep{0pt}%
    \put(0,0){\includegraphics[width=\unitlength,page=1]{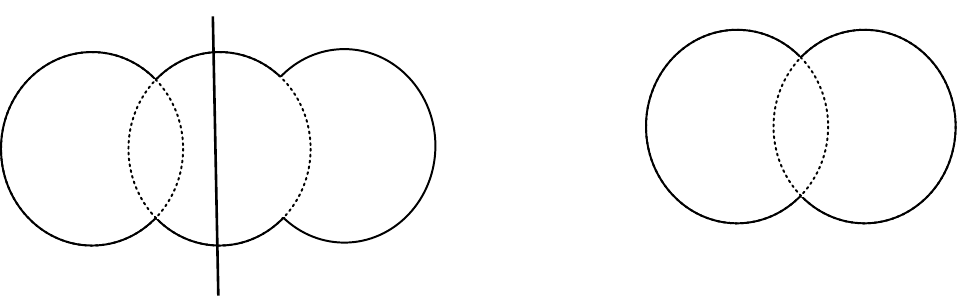}}%
    \put(0.47436169,0.1430){\color[rgb]{0,0,0}\makebox(0,0)[lt]{\lineheight{1.25}\smash{\begin{tabular}[t]{l}$\tilde\Omega_\epsilon$\end{tabular}}}}%
    \put(0,0){\includegraphics[width=\unitlength,page=2]{GPexample.pdf}}%
    \put(0.9402196,0.1430){\color[rgb]{0,0,0}\makebox(0,0)[lt]{\lineheight{1.25}\smash{\begin{tabular}[t]{l}$\Omega_\epsilon$\end{tabular}}}}%
    \put(0.19174798,0.26296637){\color[rgb]{0,0,0}\makebox(0,0)[lt]{\lineheight{1.25}\smash{\begin{tabular}[t]{l}$\tau$\end{tabular}}}}%
  \end{picture}%
\endgroup%

    \caption{$\tilde\Omega_\epsilon$ and $\Omega_\epsilon$ for $k=2$.}
    \label{gpexample}
\end{figure}
\end{proof}

We now discuss situations in which  the minimum in the right hand side of inequality \eqref{min} is achieved by $\sigma_k^N(\Om)$.\\ In \cite{BKPS}, Ba\~nuelos,  Kulczycki,  Polterovich,
and  Siudeja proved Friedlander type inequality between the eigenvalues of Steklov-Neumann and Steklov-Dirichlet eigenvalues of Euclidean domains satisfying the so-called \textit{weak} or \textit{standard John's conditions}. 

\begin{defn}\label{def.john}  
A domain $\Om$ in $\R^2$ with Lipschitz boundary is said to satisfy the \emph{standard John's condition} if  $\partial \Om\cap \{y=0\}=I$ and
 $\Om \subset I\times (-\infty,0).$\\
We say that $\Om$ satisfies the \emph{weak John's condition} if $\partial \Om\cap \{y=0\}=I$,  $\Om\subset \{(x,y)\in \R^2: y < 0\}$,
and for all $a >0$, 
\begin{equation}\label{johna}\int_\Om\, e^{ay}dy\,dx\leq \int_I\int_{-\infty}^0\,e^{ay}dy\,dx=\frac{|I|}{a}.\end{equation}

In particular, any domain that satisfies the standard John's condition also satisfies the weak John's condition.
\end{defn}

\begin{prop}\label{prop.john} Suppose $\Om\subset\R^2$ satisfies the weak John's condition. {Let $\pa\Om=\pa_S\Om \sqcup\pa_*\Om$ be a decomposition of $\pa \Om$, where $\pa_S\Om=\Om\cap\{y=0\}$.}     Then
\begin{equation}\label{johnbound}\sigma_k^N(\Om)L(\pa_S\Om) \leq (2k-1)\pi,\qquad k\ge1\end{equation}
      Moreover,  for $k=1$,  inequality \eqref{johnbound} is sharp in the sense that   
       there exists a family of domains $\Om_n$ satisfying the standard John's condition such that $\lim_{n\to\infty}\sigma_1^N(\Om_n)L(\pa_S\Om)= \pi.$
 For $k\ge2$, if $\Om$ satisfies the weak John's condition with inequality~\eqref{johna} being strict for all~$a$, then inequality~\eqref{johnbound} is strict.    
\end{prop}

\begin{proof} In \cite{BKPS}, it is proved that   that 
$\sigma_k^N(\Om)\leq\sigma_{k-1}^D(\Om)$ for any $k\ge1$.
Hence, inequality~\ref{johnbound} immediately  follows from Theorem~\ref{minpi}.\\
To show the sharpness of $k=1$, let  $
\Om_n=I\times(-n,0)$, where $I$ is an open interval. The sequence of $\Om_n$ satisfies the standard John's condition. By separation of variables (cf. \cite[\S 2.1]{BKPS}), we have
$\sigma_1^N(\Om_n)=\frac{\pi}{|I|}\tanh\left(\frac{\pi n}{|I|}\right)$. Therefore, $\sigma_1^N(\Om_n)|I|$ tends to $\pi$ as $n\to\infty$.

We now prove the last statement of the proposition.  Let  $\Om'$ be a domain satisfying same hypothesis as $\Om$ properly containing $\Om$ with $\psom=\psom'$.  Since $\sigma_k^N(\Om)<\sigma_k^N(\Om')$ (see Lemma~\ref{DoMono}), the inequality~\eqref{johnbound} must be strict for $\Om$. It completes the proof.\end{proof}

Our final example in this section addresses domains with dihedral symmetry.

\begin{ex}\label{HD remark} a) Let $\Omega$ be a compact  simply-connected Riemannian surface with Lipschitz boundary.  Assume that $\Omega$ is symmetric under an action of the dihedral group $D_4$ generated by rotation $\rho$ through angle $\frac{\pi}{2}$ about a point $p\in\Om$ and reflection $\tau$ across a geodesic $\beta$ through $p$.   With the notation in Remark \ref{nota.orb}, for ${\Om_\tau}:=\Om/\la \tau \ra$ we have 
\begin{equation}\label{ohd_max}\sigma_1^N(\Om_\tau)L(\pa_S{\Om_\tau})=\sigma_0^D(\Om) L(\psom)\leq \pi \end{equation} 
and equality holds in  (\ref{ohd_max}) only if $\Om$ is isometric to a half-disk up to conformal change of metric away from the boundary.   

Inequality~(\ref{ohd_max}) can be verified by an elementary argument following the same idea used by 
James R. Kuttler and Vincent G. Sigillito in \cite{KuSi}.  Let $\mu$ be the isometry in $D_4$ given by reflection across the geodesic through $p$ orthogonal to $\beta$.    Consider the eigenspace corresponding to the lowest non-zero Steklov eigenvalue of $\Om$.  Since each eigenspace is invariant under the action of the dihedral group and since $\tau$ and $\mu$ commute, there exists a basis of first eigenfunctions each of which is either invariant or anti-invariant with respect to $\tau$ and also either invariant or anti-invariant with respect to $\mu$.   Since a first eigenfunction must have exactly two nodal domains (as first pointed out in \cite{KuSi}), the eigenfunction must be invariant with respect to one of these reflections and anti-invariant with respect to the other.   Moreover, the rotation $\rho$ conjugates one reflection to the other.  Using the $D_4$-invariance of the eigenspace, we thus obtain two first eigenfunctions, one of which is $\tau$-invariant and the other $\tau$-anti-invariant.   These correspond to eigenfunctions for the two mixed Steklov problems on $\Om_\tau$.   We thus have $
\sigma_1^N(\Om_\tau)=\sigma_0^D(\Om_\tau)= \sigma_1(\Om).$
On the other hand $L(\pa_S\Om_\tau)=\frac{1}{2}L(\partial \Om)$. Now, inequality \eqref{ohd_max} and the equality statement follows from Weinstock's theorem.\\

\noindent b) 
We now relax the assumption of simply connectedness and assume $\Om$ has genus zero and that $\pa\Om$ has boundary decomposition of type $(b;b-1;m)$ for some $m$. We also assume $\Om$ has the same group of symmetry as in part a). Note that point $p$ may be in the ambient space and not in $\Om$. Then repeating the same line of arguments as above and using inequality \eqref{inq:gkl}, we get. 
\begin{equation}\label{ohdmax}
\sigma_1(\Om_\tau)L(\partial \Om_\tau)=\frac{1}{2} \sigma_1(\Om)L(\partial \Om)\le \frac{1}{2}\Sigma_{1}(0,2(b-1)+m).
\end{equation}

\end{ex}

\section{Eigenvalue bounds for domains with symmetry}\label{Bandle}

Although Weinstock showed that $\sigma_1(\Om)L(\pa \Om)$ is maximized when $\rho \equiv 1$ by a disk among all simply connected planar domains (see \eqref{weinstock}), the disk fails to maximize higher eigenvalues.   Instead, one has an upper bound given by the Hersch-Payne-Schiffer inequality\eqref{HPS1}, shown by Girouard and Polterovich to be sharp.  However, Bandle showed that if one restricts attention to domains with rotational symmetry of order $p$, then the disk continues to maximize $\sigma_k(\Om)L(\pa \Om)$ for all $k\leq p-1$.    
\begin{nota}\label{nota.bandle}

For any integer $p\ge2$, let $\mathcal{D}_p$ be the class of all pairs $(\Om, \rho)$ of
simply connected planar domains $\Om$ with Lipschitz  boundary and rotational symmetry of order $p$ centered at the origin, and rotationally symmetric weights $\rho$ on $\pa \Om$. For any positive integer $k$ , let 
\begin{equation}\label{sigsup}\Sigma_k^{p}=\sup_{(\Om,\rho)\in\mathcal{D}_p}\sigma_k(\Om)L_\rho(\pa\Om).\end{equation}  

\end{nota}

\begin{thm}(Bandle \cite[Theorem 3.15]{bandle})\label{thm.bandle} For any inetger $p\ge 2$
\begin{equation}\label{bandle1}
  \Sigma_k^p=\sigma_k(\disk)L(\pa \disk)=\begin{cases}  (k+1)\pi&\text{$k$ odd,\,$1\le k\le p-1$}, \\
 k\pi& \text{$k$ even,\,$2\le k\le p-1$}.
 \end{cases}  
\end{equation}
\end{thm}

Bandle assumed piecewise analytic boundary, but it suffices to assume that the boundary is  Lipschitz.   Indeed, every Lipschitz domain $\Om$ can be approximated by a sequence of domains $\{\Om_n\}{\subset\mathcal{D}_p}$ with piecewise analytic boundary satisfying the assumption of Proposition \ref{Proposition 2.3, BuNa2.3}: take $\Om_n$ to have polygonal boundary with vertices lying on $\pa\Om$ such that the distance between vertices goes to zero as $n$ goes to infinity. The rotationally symmetric $\rho$ can also be estimated by a sequence of smooth rotationally symmetric functions so, without loss of generality, we assume $\rho\equiv1$ unless otherwise stated.  Note that the case $k=p-1$ in  inequality \eqref{bandle1} was first proven in \cite{HPS}.  

\label{Akhmetgaliyev} In \cite{AKO}, Akhmetgaliyev, Kao and Osting carried out numerical optimization of the Steklov eigenvalues of planar domains normalized by area, leading them to conjecture that the area-constrained maximizer $\Om_p^*$ of the $p$th eigenvalue is unique and has $p$-fold rotational symmetry and at least one axis of reflection symmetry.  (See also \cite{Bo17}.)  For each $p\leq 100$, they obtained numerical estimates of the extreme value of $\sigma_p$ with area constrained to be one and also estimates of several higher order eigenvalues of the numerically maximizing domain $\Om_p^*$.  In Table \ref{table1}, we have used the isoperimetric inequality, $\sigma_k(\Om_p^*)L(\pa\Om_p^*)>2\sqrt{\pi}\sigma_k(\Om_p^*){\rm Area}(\Om_p^*)^{1/2}$, to convert their estimates in \cite[Fig 2]{AKO} to lower bounds for  $\sigma_k(\Om_p^*)L(\pa \Om_p^*)$.   Letting 
 \[K(p):=\{k\in \Z^+: \Sigma^p_k=\sigma_k(\disk )L(\partial \disk )\},\] 
 these estimates show that $p, p+1\notin K(p)$ when  $p\le 10$.
    \begin{center}
    \begin{table}
    \begin{tabular}{c|cccccccc}
        $k/p$ & $3$&$4$&$5$&$6$&$7$&$8$&$9$&$10$ \\ \hline
         $3$&${\bf 4.67\pi}$&&&&&&& \\
         $4$ &${\bf 4.67\pi}$&${\bf5.96\pi}$&&&&&& \\
          $5$ &$4.67\pi$&${5.96\pi}$&${\bf7.33\pi}$&&&&& \\
           $6$ &$5.54\pi$&${\bf 6.14\pi}$&${\bf7.33\pi}$&${\bf8.62\pi}$&&&& \\
            $7$ &$6.79\pi$&$6.14\pi$&${7.33\pi}$&${\bf8.62\pi}$&${\bf9.98\pi}$&&& \\
            $8$ &$6.79\pi$&$7.32\pi$&$7.59\pi$&${\bf 8.77\pi}$&${\bf9.98\pi}$&${\bf11.28\pi}$& \\
            $9$ &$8.6\pi$&$8.27\pi$&$7.59\pi$&$8.77\pi$&${9.98\pi}$&${\bf11.28\pi}$&${\bf12.63\pi}$& \\
            $10$ &$8.6\pi$&$8.27\pi$&$9.17\pi$&$9\pi$&${\bf 10.21\pi}$& ${\bf 11.4\pi}$&${\bf12.63\pi}$&${\bf13.94\pi}$ \\
             $11$ &$10.10\pi$&$9.74\pi$&${9.93\pi}$&$9\pi$&$10.21\pi$&$11.4\pi$&${\bf12.63\pi}$&${\bf13.94\pi}$ \\
              $12$ &$10.31\pi$&${10.19\pi}$&${9.93\pi}$&${11.65\pi}$&${10.98\pi}$&${11.65\pi}$&${\bf12.82\pi}$&{$\bf14.03\pi$} \\
    \end{tabular} 
    \bigskip
 
    \caption{Lower bounds of $\sigma_k(\Om_p^*)L(\pa\Om_p^*)$ for $3\le p\le 10,~ 3\le k\le12$. The numbers in bold are strictly greater than $\sigma_k(\disk)L(\partial \disk)$. 
    }\label{table1}
    \end{table} 
\end{center}

Contrasting Bandle's Theorem \ref{thm.bandle} with the Hersch-Payne-Schiffer inequality, we ask:

\begin{ques}\label{sharp?} Is Bandle's Theorem sharp; i.e., for fixed $p$, is $p-1$ the maximum value of $k \in K(p)$?
\end{ques}

\begin{ques}\label{hps?}  Does  there exist $k_0\geq p$ such that the Hersch-Payne-Schiffer  inequality is sharp on $\mathcal{D}_p$ for $k\geq k_0$, i.e. is it the case that
$\Sigma_k^p=2\pi k$ once $k\ge k_0$ is sufficiently large? If so, is $k_0=p$ big enough?
 
\end{ques}

When $p=2$, Girouard and Polterovich's proof \cite{GP} that the Hersch-Payne-Schiffer inequality is sharp for every $k\geq 2$ yields affirmative answers to both Question~\ref{sharp?} and Question~\ref{hps?}, with $k_0=p=2$ in the latter.  In particular, they constructed for each $k\geq 2$ a family of domains $\{\Om_{k,\epsilon}\}$ such that $\lim_{\epsilon\to 0}\,\sigma_k(\Om_{k,\epsilon})=2\pi k$.   Their domains have rotational symmetry of order 2.  Thus
   
$$\Sigma^2_k=2\pi k \mbox{ for all }k\geq 2.$$   

This construction motivates the proof of the next theorem, which gives lower bounds on $\Sigma_k^p$ that we will use to address the questions above. We postpone the proof until after the statements of some consequences of the theorem.

\begin{thm}\label{bndlesharp} In the notation of~\ref{nota.bandle} and Equation~(\ref{sigsup}), we have 
\begin{enumerate}[(i)]
\item $\Sigma_k^p> p \pi$, for all $k\geq p$.

\item $\Sigma_k^p\geq (2\left\lfloor\frac{k}{p}\right\rfloor -1) p\pi$, for all $k\geq 2p$. Here $\lfloor \cdot\rfloor$ is the greatest integer function.

\item When $3\leq p\leq 6$, we have $\Sigma_k^p= 2k\pi$ when $k>p$ satisfies $k\equiv 1$ mod $p$.   
\end{enumerate} 
\end{thm}

The lower bounds in Theorem \ref{bndlesharp} are not strong enough to answer Question~\ref{hps?}, but part (ii) leads to the following weaker statement:

\begin{cor}\label{BandleLimit}
For every $p$, we have 
$$\lim_{k\to\infty}\,\frac{\Sigma^p_k}{2\pi k}=1.$$
\end{cor}

 Theorem \ref{bndlesharp} yields more substantial results concerning Question~\ref{sharp?}: 

  \begin{cor}\label{cor3.3.3} Recall that
 $K(p)=\{k\in \Z^+: \Sigma^p_k=\sigma_k(\disk)L(\partial \disk)\}.$
We have:
 \begin{enumerate}[(i)]
 \item $p\notin K(p)$, when $p$ is even.
  \item $K(p)\cap [p+1,\infty)=\emptyset$, when $2\leq p\leq 6$.
 \item $K(p)\cap [2p,\infty)=\emptyset$, when $p$ is odd.
 \item $K(p)\cap [2p,\infty) \subset \{3p-1\}$, when $p$ is even.
 \end{enumerate}
  \end{cor}

 \begin{proof}  Items (i) and (ii) follow from Theorem~\ref{bndlesharp}, items (i) and (iii), respectively.
 
 We prove items (iii) and (iv) jointly.  
 For arbitrary $p$ and for $k\geq 2p$, write $k=mp+q$ with $1\leq q\leq p-1$.   Theorem~\ref{bndlesharp}(ii) says that 
 \begin{equation}\label{kpeq}\Sigma^p_k\geq \big(k+(m-2)p +(p-q)\big)\pi\geq 2\pi\left\lfloor \frac{k+1}{2}\right\rfloor =\sigma_k(\disk)L(\pa \disk).\end{equation} 
  Since $p-q\geq 1$ and $m\geq 2$, the second inequality in~\eqref{kpeq} is strict except possibly when $m=2$ and $p-q=1$, i.e., when $k=3p-1$ with $p$ even.   
 \end{proof}
 
 It seems highly likely that $K(p)\cap [2p,\infty) =\emptyset$ for all even as well as odd $p$.   Indeed,  if $3p-1$ actually lies in $K(p)$, then we have $\Sigma^p_{3p-1}=3\pi p$, while $\Sigma^p_{2p}\geq 3\pi p$ by Theorem~\ref{bndlesharp}.  Thus in this case, we must have $\Sigma^p_{2p}=\Sigma^p_{2p+1}=\dots =\Sigma^p_{3p-1}$.  Moreover, we must also have $p\geq 8$ (by item (ii) of Corollary \ref{cor3.3.3}), so we have a sequence of at least eight values of $k$ on which the map $k\mapsto \Sigma^p_k$ remains constant, certainly an unlikely scenario! We conjecture that \\
 
  \noindent\textbf{Conjecture.~} $K(p)\cap[p,\infty)=\emptyset$.

  \begin{proof}[Proof of Theorem~\ref{bndlesharp}]
(i) Let $k=p$. Let  $G_p$ be a regular $p$-gon centered at the origin with sides of length $2$ and with vertices lying on the rays given in polar coordinates by $\left\{(\theta, r):\theta=\frac{2\pi i}{p}\right\}_{i=1}^p$. To each side, we attach a unit half disk.  We denote  the resulting domain  by $\Om$. See Figure \ref{threeballs}~(a). 
\begin{figure}
    \centering
  \hspace{-1cm}  
\begingroup%
  \makeatletter%
  \providecommand\color[2][]{%
    \errmessage{(Inkscape) Color is used for the text in Inkscape, but the package 'color.sty' is not loaded}%
    \renewcommand\color[2][]{}%
  }%
  \providecommand\transparent[1]{%
    \errmessage{(Inkscape) Transparency is used (non-zero) for the text in Inkscape, but the package 'transparent.sty' is not loaded}%
    \renewcommand\transparent[1]{}%
  }%
  \providecommand\rotatebox[2]{#2}%
  \newcommand*\fsize{\dimexpr\f@size pt\relax}%
  \newcommand*\lineheight[1]{\fontsize{\fsize}{#1\fsize}\selectfont}%
  \ifx\svgwidth\undefined%
    \setlength{\unitlength}{258.74347897bp}%
    \ifx\svgscale\undefined%
      \relax%
    \else%
      \setlength{\unitlength}{\unitlength * \real{\svgscale}}%
    \fi%
  \else%
    \setlength{\unitlength}{\svgwidth}%
  \fi%
  \global\let\svgwidth\undefined%
  \global\let\svgscale\undefined%
  \makeatother%
  \begin{picture}(1,0.70546415)%
    \lineheight{1}%
    \setlength\tabcolsep{0pt}%
    \put(0,0){\includegraphics[width=\unitlength,page=1]{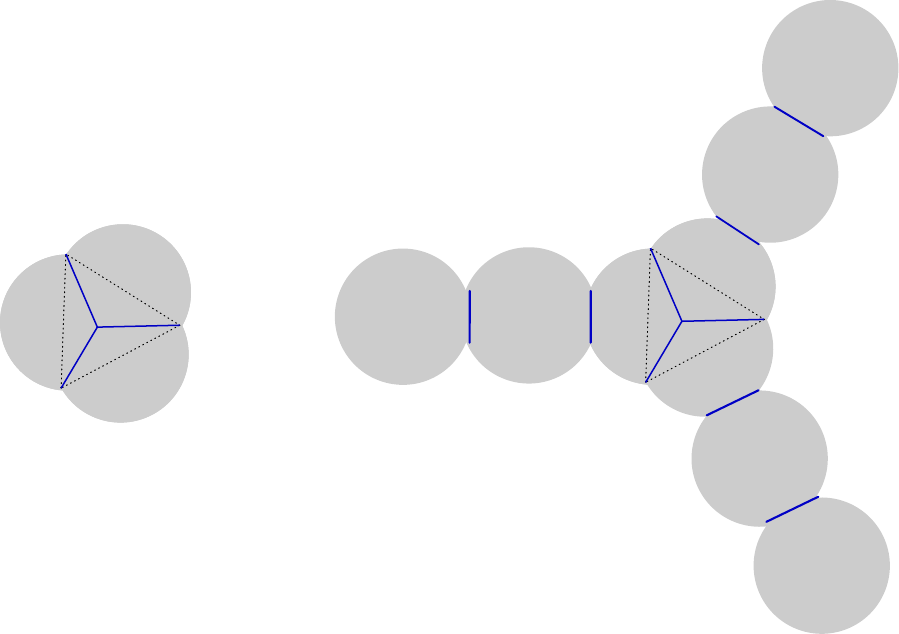}}%
    \put(0.094,0.15){\color[rgb]{0,0,0}\makebox(0,0)[lt]{\lineheight{1.25}\smash{\begin{tabular}[t]{l}$(a)$\end{tabular}}}}%
    \put(0.57193061,0.15){\color[rgb]{0,0,0}\makebox(0,0)[lt]{\lineheight{1.25}\smash{\begin{tabular}[t]{l}$(b)$\end{tabular}}}}%
  \end{picture}%
\endgroup%

   \caption{Picture $(a)$ shows domain $\Om$ when $p=3$. Picture $(b)$  shows $\Om_\epsilon$ when $p=3$ and $k=9$. We impose the Neumann boundary condition on blue lines.}
    \label{threeballs}
\end{figure}
Consider the sectors $C_p^i := \left\{(\theta, r): \frac{2\pi (i-1)}{p}\le \theta\le\frac{2\pi i}{p} \right\}$ in the plane. Define \label{omi}$\Om_i:=\Om\cap C_p^i$ and consider the Steklov-Neumann problem on $\Om_i$ with Neumann condition on $\Om\cap \pa C_p^i$ and Steklov condition on the half-circle. By Lemma \ref{DoMono}, we have 
 \begin{equation}\label{smonoton}\sigma_p(\Om)\ge\min\{\sigma^N_1(\Om_1),\cdots,\sigma^N_1(\Om_p)\}. \end{equation}

Since the half-disk $\hD(\pi)$  is a proper subset of $\Om_i$, we can impose the Neumann condition on the diameter of $\hD(\pi)$ and use the monotonicity of Steklov-Neumann eigenvalues (see Lemma~\ref{DoMono}) to conclude that
 \begin{equation}\label{ahd+}\sigma_1^N(\Om_i)>\sigma_1^N(\hD)=1.\end{equation}
The $\Om_i$s are mutually isometric, so
 $\sigma_p(\Om)L(\pa\Om)>\pi p$,
proving (i).
 
(ii) Let $k=mp$ with $m\ge 2$. Let \[U_{\epsilon}:=\bigcup_{j=1}^{m-1}D_{\epsilon,j},\]
where $D_{\epsilon,j}$ is a disk of radius $1+ \epsilon$ centered at the point $(2j,0)$.    Let $G_p$ be the polygon defined in the proof of (i), and let $l$ be the distance of each side of $G_p$
from the origin.   Let $U_\epsilon^i$, $i=1,\ldots p$ be the image of $U_\epsilon$ under the translation  $(x,y)\mapsto(x+l,y)$ and  the rotation  through angle $\frac{\pi {(2i-1)}}{p}$.   We define $$ \Om_\epsilon:=\bigcup_i^{p}U_\epsilon^i\cup\Om,$$
where $\Om$ is the domain constructed in the proof of (i).  See Figure \ref{threeballs}~(b) for an illustration.

We now partition $\Om_\epsilon$ into $mp$ domains $\{\Om_{\epsilon,j}\}_{j=1}^{mp}$ where each domain is isometric to one of the following  domains (see Figure \ref{vdomains}):
 \begin{itemize}
 \item $V_{1,\epsilon}=\tilde\Om\cap\{x\le l+1-\epsilon\}\cap \{\frac{-\pi}{p}\le \theta\le\frac{\pi}{p} \}$, where $\tilde\Om$ is obtained by rotating $\Om$ by $\frac{\pi}{p}$ around the origin   (here $(x,y)$ are the usual Cartesian coordinates);
  \item $V_{2,\epsilon}=D_{\epsilon,1}\cap\{1\le x\le 3\};$
 \item $V_{3,\epsilon}:=D_{\epsilon,1}\cap\{x\ge 1\}.$
 \end{itemize}

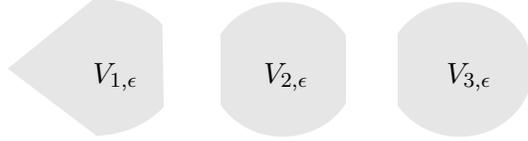
\begin{figure}
    \centering
  \hspace{-1cm}  
\begingroup%
  \makeatletter%
  \providecommand\color[2][]{%
    \errmessage{(Inkscape) Color is used for the text in Inkscape, but the package 'color.sty' is not loaded}%
    \renewcommand\color[2][]{}%
  }%
  \providecommand\transparent[1]{%
    \errmessage{(Inkscape) Transparency is used (non-zero) for the text in Inkscape, but the package 'transparent.sty' is not loaded}%
    \renewcommand\transparent[1]{}%
  }%
  \providecommand\rotatebox[2]{#2}%
  \newcommand*\fsize{\dimexpr\f@size pt\relax}%
  \newcommand*\lineheight[1]{\fontsize{\fsize}{#1\fsize}\selectfont}%
  \ifx\svgwidth\undefined%
    \setlength{\unitlength}{214.33801329bp}%
    \ifx\svgscale\undefined%
      \relax%
    \else%
      \setlength{\unitlength}{\unitlength * \real{\svgscale}}%
    \fi%
  \else%
    \setlength{\unitlength}{\svgwidth}%
  \fi%
  \global\let\svgwidth\undefined%
  \global\let\svgscale\undefined%
  \makeatother%
  \begin{picture}(1,0.24366993)%
    \lineheight{1}%
    \setlength\tabcolsep{0pt}%
    \put(0,0){\includegraphics[width=\unitlength,page=1]{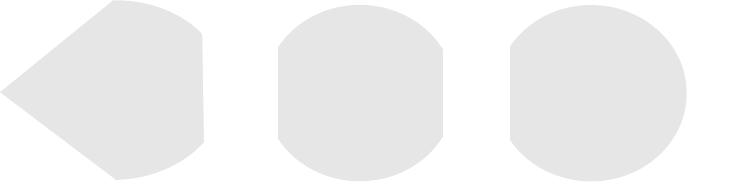}}%
        \put(0.1498461,0.095252){\color[rgb]{0,0,0}\makebox(0,0)[lt]{\lineheight{1.25}\smash{\begin{tabular}[t]{l}$V_{1,\epsilon}$\end{tabular}}}}%
    \put(0.4498461,0.095252){\color[rgb]{0,0,0}\makebox(0,0)[lt]{\lineheight{1.25}\smash{\begin{tabular}[t]{l}$V_{2,\epsilon}$\end{tabular}}}}%
    \put(0.77164048,0.095368619){\color[rgb]{0,0,0}\makebox(0,0)[lt]{\lineheight{1.25}\smash{\begin{tabular}[t]{l}$V_{3,\epsilon}$\end{tabular}}}}%
  \end{picture}%
\endgroup%

  \caption{Domains $V_{j,\epsilon}$.}
    \label{vdomains}
\end{figure}
We impose the Neumann boundary condition on the flat parts of  $\partial V_{i,\epsilon}$. 
Again, by Lemma \ref{DoMono}, we have 
\[\sigma_{mp}(\Om_\epsilon)\ge\min\{\sigma_1^N(\Om_{\epsilon,j}):j=1,\ldots,mp\}=\min\{\sigma_1^N(V_{1,\epsilon}),\sigma_1^N(V_{2,\epsilon}),\sigma_1^N(V_{3,\epsilon})\}.\] 
Using Proposition \ref{Proposition 2.3, BuNa2.3} and inequality \eqref{ahd+}, we have $$\lim_{\epsilon\to 0}\sigma_1^N(V_{1,\epsilon})=\sigma_1^N(\Om_1)>1,$$ where $\Om_1$ is defined as in the proof of 
case (i) and $$\lim_{\epsilon\to0} \sigma_1^N(V_{2,\epsilon})=\lim_{\epsilon\to0}\sigma_1^N(V_{3,\epsilon})=\sigma_1(D)=1.$$
We complete the proof of (ii) by noting that $$\lim_{\epsilon\to 0}L(\pa \Om_{\epsilon})= (m-1)2\pi p+p\pi=(2m-1)\pi p,$$
giving
$$\lim_{\epsilon\to 0}\sigma_{mp}(\Om_\epsilon) L(\pa \Om_{\epsilon})\geq (2m-1)\pi p.$$

(iii) Let $k=mp+1$, where $m\ge1$ and $2\le   p\le 6$.  Since the Hersch-Payne-Schiffer inequality implies that $\Sigma_k\leq 2\pi k$, it suffices to prove the opposite inequality.

Let $D_\epsilon$ be the disk of radius $1+2\epsilon$ centered at the origin and let
\[U_\epsilon:=\bigcup_{j=1}^{m}D_{j,\epsilon}\]
where $D_{j,\epsilon}$ is the disk of radius $1+\epsilon$ centered at $(2j+\epsilon,0)$.  Let $U_\epsilon^i$, $i=1,\ldots, p$ be the image of $U_\epsilon$ under  the rotation  by $\frac{(2\pi (i-1))}{p}$.  Since $2\le p\le 6$, and since $D_\epsilon$ has radius greater than the radii of $U^i_\epsilon$, the $U^i_\epsilon$ are mutually disjoint for $\epsilon$ small enough. We define $\Om_\epsilon:= (\cup_{i}U^i_\epsilon)\cup D_\epsilon$. See Figure \ref{fig:om_e/o}.
\begin{figure}
    \centering
    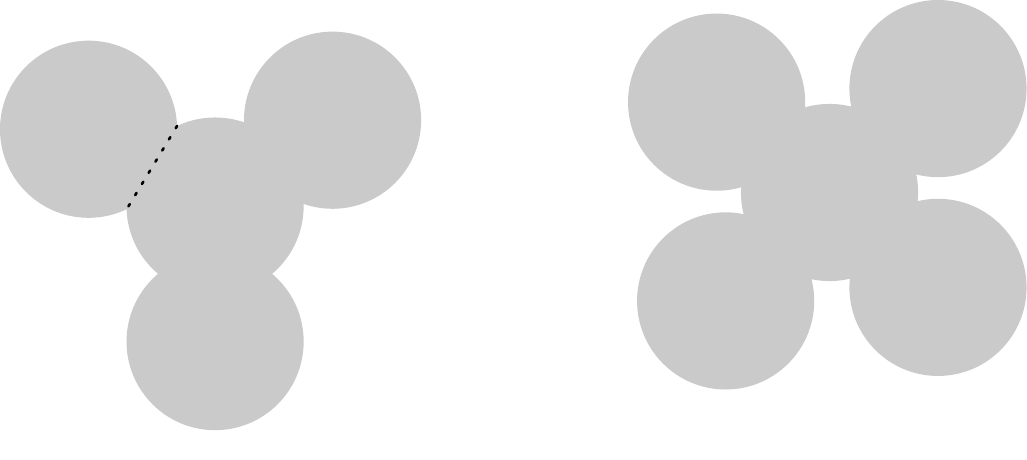
    \caption{In these two examples, we have $\lim_{\epsilon\to0}\sigma_4(\Om_\epsilon)L(\pa\Om_\epsilon)=8\pi$ for the picture on the left and $\lim_{\epsilon\to0}\sigma_5(\Om_\epsilon)L(\pa\Om_\epsilon)= 10\pi$ for the picture on the right. Let $\Om_{\epsilon,\tau}=\langle\tau\rangle\backslash\Om_\epsilon$. We have $\lim_{\epsilon\to0}\sigma_3(\Om_{\epsilon,\tau})L(\pa\Om_{\epsilon,\tau})=4\pi$ for the picture on the left and $\lim_{\epsilon\to0}\sigma_3(\Om_{\epsilon,\tau})L(\pa\Om_{\epsilon,\tau})= 5\pi$ for the picture on the right.}
    \label{fig:om_e/o}
\end{figure}
We partition $\Om_\epsilon$ into $mp+1$ domains $\{\Om_{\epsilon,j}\}_{j=1}^{mp+1}$ by inserting line segments joining the two points of intersection of each pair of overlapping disks.   Thus each $\Om_{\epsilon,j}$ is a disk with one or more arcs of the boundary circle replaced by chords.  We impose the Neumann boundary condition on the chords and Steklov conditions on the circular parts of the boundaries. Note that by \eqref{buna} $$\lim_{\epsilon \to 0}\sigma^N_1(\Om_{j,\epsilon})=\sigma_1(\disk)=1,\quad j=1,\ldots,mp+1.$$  Thus $\displaystyle{\lim_{\epsilon \to 0}} \ \sigma_{mp+1}(\Om_\epsilon)\geq 1$  while $\displaystyle{\lim_{\epsilon\to0}} \ L(\pa\Om_\epsilon)=2\pi(mp+1),$
giving  $$\lim_{\epsilon\to0}\sigma_{k}(\Om_{\epsilon})L(\pa\Om_{\epsilon})\ge 2\pi(mp+1).$$  
This completes the proof of (iii).   

 Notice that this construction of $\Om_\epsilon$ does not work for $p\ge 7$ because the $U_\epsilon^i$ will not be disjoint and $\Om_\epsilon$ may not remain simply connected. 
 \end{proof}
 
We now extend Bandle's Theorem \ref{thm.bandle} to domains whose doubles are in $\mathcal{D}_p$ and consider the Steklov-Neumann problem on these domains.  
 
\begin{nota}\label{nota.bandle2}
 For any integer $p\ge2$, let $\mathfrak{D}_p$ be the class of all planar Lipschitz domains $\Om$ whose boundary is of the following form:  $\partial\Om=\psom \cup\pnom$ where
$\pnom$ is an interval symmetric about zero on the $x$-axis, and 
    the double of $\Om$ over $\pnom$ is a planar domain lying  in $\mathcal{D}_p$. 
 For $k \in \Z^+$, let 
\begin{equation}\label{prefl}\mathfrak{S}_k^p:=\sup_{\Om\in\mathfrak{D}_p}\sigma_k^N(\Om)L(\pa_S\Om).\end{equation}
 \end{nota}
 
We can view $\Om$ as the quotient of its double   $\tilde\Om\in \mathcal{D}_p$ by its reflection $\tau$ in $\pa_N\Om$. By Remark \ref{nota.orb}, $\sigma_k^N(\Om)$ is equal to the Steklov eigenvalues of $\tilde\Om/\langle\tau\rangle$ as an orbifold. 
 \begin{cor}\label{wimproved}  In the notation of \ref{nota.bandle2}, for every $\Om\in \mathfrak{D}_p$ and every $0\le k< \lfloor\frac{p+1}{2}\rfloor$, we have \begin{equation}\label{imbandle2}
 \mathfrak{S}_k^p=\sigma_k^N(
    \hD)L(\pa \hD)= k \pi.
 \end{equation}
 \end{cor}

 \begin{proof}
With a similar argument as in the proof of Proposition \ref{minpisharp}, we have $\sigma_k(\Om)\le \sigma_{2k}(\tilde\Om)$. Hence, by Bandle's Theorem \ref{thm.bandle} we conclude
 \[\sigma_k(\Om)L(\pa\Om)\le \sigma_{2k}(\tilde\Om)\frac{L(\pa\tilde\Om)}{2}\le k\pi,\qquad 0\le 2k\le p-1.\]
  \end{proof}
We can now ask Questions \ref{sharp?} and \ref{hps?}, replacing $\Sigma_k^p$ by $\mathfrak{S}_k^p$. Theorem \ref{dpmodreflection} is the counterpart of Theorem \ref{bndlesharp}. It leads to Corollaries \ref{BandleLimit} and  \ref{cor3.3.8}, which give partial answers to those questions. 
  \begin{thm}\label{dpmodreflection} Let $\mathfrak{p} :={\lfloor\frac{p+1}{2}\rfloor}$. In the notation of \ref{nota.bandle2} we have 
  \begin{enumerate}[(i)]
  \item {$\mathfrak{S}^p_k>\frac{p\pi}{2}$, for all $k\ge\mathfrak{p}.$}
  \item $\mathfrak{S}^p_k\ge\left(\left\lfloor\frac{k}{\mathfrak{p}}\right\rfloor-\frac{1}{2}\right)\pi p,$ for every $k\ge2\mathfrak{p}$.  
  \item When $3\le p\le 6$, we have $\mathfrak{S}^p_k= \left(\left\lfloor\frac{k}{\mathfrak{p}}\right\rfloor p+1\right)\pi,$ for every $k\ge 2\mathfrak{p}$ such that $k\equiv 1$ mod $\mathfrak{p}$. 
  \item When $p=2$,  $\mathfrak{S}^p_{k}= 2k\pi,$ for every $k\ge1$.
  \end{enumerate}
  \end{thm}
As an immediate consequence of part (ii) and the Hersch-Payne-Schiffer inequality, we have
\begin{cor}\label{BandleLimit2}
For every $p$, we have 
$$\lim_{k\to\infty}\,\frac{\mathfrak{S}^p_k}{2\pi k}=1.$$
\end{cor}

From Theorems \ref{bndlesharp} and \ref{dpmodreflection}, we obtain
 \begin{cor}\label{cor3.3.8} Let
 $\mathcal{K}(p):=\{k\in \Z^+: \mathfrak{S}^p_k=\sigma_k(\hD(\pi))L(\partial \hD(\pi))=\pi k\}. $ Again using the notation $\mathfrak{p}=\left\lfloor\frac{p+1}{2}\right\rfloor$, we
 then have
 \begin{enumerate}[(i)]
     \item $\frac{p}{2}\notin\mathcal{K}(p)$, when $p$ is even.

     \item $\mathcal{K}(p)\cap [p,\infty)=\emptyset$, when $p$ is even.
     \item $\mathcal{K}(p)\cap [3\mathfrak{p},\infty)=\emptyset$, when $p\ge3$ is odd.
     \item $\mathcal{K}(p)\cap (\mathfrak{p}, \infty)=\emptyset$, when $2\le p\le 6$.
 \end{enumerate}
 \end{cor}
 
\noindent\textbf{Conjecture.~} $\mathcal{K}(p)\cap[\mathfrak{p},\infty)=\emptyset$, where $\mathfrak{p}=\left\lfloor\frac{p+1}{2}\right\rfloor$.
  
    \begin{proof}[Proof of Theorem \ref{dpmodreflection}]
For Parts (i)-(iii), we consider the domains used in the proof of Theorem \ref{bndlesharp}. All have a reflection symmetry $\tau$ passing through the origin and without loss of generality we can assume the axis of symmetry is the $x$-axis. The partition given for the domain $\tilde \Om$ in the proof of Theorem \ref{bndlesharp} gives rise to a partition for $\Om$ when we view it as $\tilde\Om_{\tau} := \tilde\Om/\langle \tau \rangle$. 
For the reader's convenience we summarise the proof of each part. 

(i) Let $\{\tilde\Om_i\}$ be the partition of $\tilde \Om$ given in the proof of Theorem \ref{bndlesharp}. Then $\{\tilde\Om_{i,\tau}:=\tilde\Om_i/\langle\tau\rangle\}$ gives a partition for $\tilde \Om_\tau$ which, for a given $p$, contains $\mathfrak{p}$ elements.   Each element is either $\tilde\Omega_i$ defined on page \pageref{omi} (if does not intersect the axis of symmetry) or $ \tilde\Omega_{i,\tau}$ when  $\tau$ fixes the axis of symmetry of $\tilde\Omega_i$. We view each set as a subset of $\Om$ and assume the Neumann boundary condition on the edge fixed by $\tau$. By the monotonicity argument  we have $$\sigma_1^N(\tilde\Omega_i)>\sigma_1(\hD)=1,\qquad \sigma_1^N(\tilde\Omega_{i,\tau})>\sigma_1(\qD)=2>1.$$
By Lemma \ref{DoMono}, $\sigma_{\mathfrak{p}}^N(\Om)$ is bounded below by $\min \left\{\sigma_1^N(\tilde\Omega_i), \sigma_1^N(\tilde\Omega_{i,\tau})\right\}=1.$  Thus
$\sigma_{\mathfrak{p}}^N(\Om)L(\pa_S\Om)>\frac{\pi p}{2}.$

(ii) For simplicity of argument, we consider odd and even $p$ separately. If $p$ is even and $k=\frac{mp}{2}$, then the partition of $\Om_\epsilon$ into $mp$ elements as in Figure \ref{vdomains} give rise to a partition for $\Om_{\epsilon,\tau}=\Om_\epsilon/\langle\tau\rangle$ with $\frac{mp}{2}$ elements. Each element in the partition is isometric to $V_{i,\epsilon}$ or $V_{i,\epsilon}/\langle\tau\rangle$, $i=1,2,3$; the latter occurs when the axis of symmetry of  $V_{i,\epsilon}$ is fixed by $\tau$. We identify a quotient space with its underlying topological space and put the Neumann boundary condition on the reflector edge; the quotients belong to $\mathfrak{D}_p$.  We have that $\sigma_k^N(\Om_{\epsilon,\tau})$ is bounded below by $$\min\{ \sigma_1^N(V_{i,\epsilon}),\sigma_1^N(  V_{i,\epsilon}/\langle\tau\rangle)| i=1,2,3\}.$$ 
The lower bound for  $\sigma_1^N(V_{i,\epsilon})$ is shown to be 1 in the proof of Theorem \ref{bndlesharp}, and 
$\lim_{\epsilon\to0}\sigma_1^N( V_{i,\epsilon}/\langle\tau\rangle)$ is bounded below by $\min\{\sigma_{1}^N(\Omega_{i,\tau}), \sigma_1(\hD(\pi))\}\ge1$.  We conclude that 
\[\lim_{\epsilon\to0}\sigma_k(\Om_{\epsilon,\tau})L(\pa_S\Om_{\epsilon,\tau})\ge \left(\frac{2k}{p}-\frac{1}{2}\right)\pi p.\]
If $p$ is odd and $k=m(\frac{p+1}{2})$, then a similar partition of $\Om_\epsilon$ into $mp$ elements gives rise to a partition $\Om_{\epsilon,\tau}$ with $m(\frac{p+1}{2})$ elements. Repeating the same argument as above, we get 
\[\lim_{\epsilon\to0}\sigma_k(\Om_{\epsilon,\tau})L(\pa_S\Om_{\epsilon,\tau})\ge  \left(\frac{2k}{p+1}-\frac{1}{2}\right)\pi p.\]
We conclude that for any $k\ge 2\mathfrak{p}$
\[\lim_{\epsilon\to0}\sigma_k(\Om_{\epsilon,\tau})L(\pa_S\Om_{\epsilon,\tau})\ge  \left(\left\lfloor\frac{k}{\mathfrak{p}}\right\rfloor-\frac{1}{2}\right)\pi p.\]
(iii) We consider the domain $\Omega_\epsilon$ and its partition into $mp+1$ domains as described in Theorem \ref{bndlesharp}(iii). It gives a partition of $\Omega_{\epsilon,\tau}$ into $m\mathfrak{p}+1$ domains. Hence for $k=m\mathfrak{p}+1$, we have $\displaystyle{\lim_{\epsilon\to0}}\sigma_{k}(\Om_{\epsilon,\tau})L(\pa_S\Om_{\epsilon,\tau})=(mp+1)\pi$ and the statement follows.\\ 
\noindent (iv)  For $p=2$, we consider the domain $\Om_\epsilon$ defined in \cite{GP09}: $\Om_\epsilon=\bigcup_{j=-m+1}^{m}D_{\epsilon,j},$
  where $D_{\epsilon,j}$ is a disk of radius $1+ \epsilon$ centered at the point $(2j-1,0)$. The partition is $\Om_j=\Om_\epsilon\cap\{2j\le x\le 2(j+1)\}$ and the symmetry $\tau$ is with respect to the $y$-axis.  We conclude that  $\displaystyle{\lim_{\epsilon\to0}}\sigma_{k}(\Om_{\epsilon,\tau})L(\pa_S\Om_{\epsilon,\tau})=2k\pi$. In order to view the domains as elements of $\mathfrak{D}_p$, we can rotate them by 90 degrees.
  \end{proof}
  
\section{Asymptotics}\label{sec:4}

We remind the reader of Notation \ref{nota:basic} and Model Example \ref{models}.

\begin{nota}\label{nota.disk vec}~
(i) Given a multiset $\mL =(\ell_1,\dots, \ell_r)$ where each $\ell_j\in \R^+$, let 
$$\disk(\mL)=\bigsqcup_{j=1}^r\,\disk(\ell_j),\,\,\,\hD(\mL)=\bigsqcup_{j=1}^r\,\hD(\ell_j),\,\mbox{and}\,\qD(\mL)=\bigsqcup_{j=1}^r\,\qD(\ell_j)$$
 denote disjoint unions of disks, half-disks or quarter disks, respectively.

(ii) Given a multiset $\mL = (\ell_1,\dots, \ell_r)$, we denote by $2\mL$ the multiset
$$2\mL=(2\ell_1,\dots, 2\ell_r).$$

(ii) Given sequences $A:=\{a_j\}_{j=0}^\infty$ and $B:=\{b_j\}_{j=0}^\infty$, we write $A\sim B$ to mean $a_j-b_j=O(j^{-\infty})$.
\end{nota}

\begin{nota}\label{nota.type}  Let $(\Om,g)$ be a connected, compact Riemannian surface with Lipschitz boundary and boundary decomposition $\partial\Om=\psom\sqcup \pnom\sqcup \pdom$ as in~\ref{nota:basic}.   
Each component of $\pa_S\Om$ that is a topological interval has one of the  following types: 
\begin{enumerate}
     \item[(D)] Both endpoints meet $\pa_D\Om$.
     \item[(N)] Both endpoints meet $\pa_N\Om$.
      \item[(DN)] One of the endpoints meets $\pa_D\Om$ and the other  $\pa_N\Om$. 
\end{enumerate}
Let $\ell_1,\dots, \ell_n$ be the lengths of the boundary components (repeated by multiplicity) and set $\mL=(\ell_1,\dots, \ell_n)$.  We will denote by $\mL_D$, $\mL_N$, and $\mL_{DN}$ the multisets of lengths of the components of $\psom$ that are topological intervals of each type.   E.g., if there are $q$ components of type $DN$ with lengths $\ell_1,\dots,\ell_q$, then 
$\mL_{DN}=(\ell_1,\dots, \ell_q)$.    
Similarly we denote by $\mL_S$ the multiset of lengths of all components of $\psom$ that are topological circles.

We will refer to the collection $\mL_S, \mL_D,\mL_N,\mL_{DN}$ as the \emph{boundary data} for $(\om,g)$ with the given boundary decomposition.   

\end{nota}

Girouard, Parnovski, Polterovich and Sher \cite{GPPS.14} obtained the asymptotics of the Steklov spectrum for Riemannian surfaces $(\om,g)$ with smooth boundary: 
\begin{equation}\label{eq:GPPS} \Stek(\om,g)\sim \Stek(\disk(\mL)).
\end{equation}
Motivated by these asymptotics and by Theorem 5.2 of \cite{ADGHRS}, we will show:  

\begin{thm}\label{asymptotic-dir-Neu} Let $(\Om,g)$ be a connected, compact Riemannian surface with Lipschitz boundary and boundary decomposition $\partial\Om=\psom\sqcup \pnom\sqcup \pdom$.    (We allow the case that $\pnom$ or $\pdom$ is empty.)  In addition to the standing assumption in~\ref{nota:basic}, assume:
\begin{itemize}
\item $\psom$ is smooth;
\item at any common endpoint $p$ of $\psom$ with either of $\pnom$ or $\pdom$, a small segment of $\pnom$, respectively $\pdom$, emanating from $p$ is a geodesic segment orthogonal to $\psom$ at $p$. 
\end{itemize} 
 Then in the notation of~\ref{nota:basic}, \ref{nota.disk vec} and \ref{nota.type}, we have 
\begin{equation}\label{thm:mix asymptotic}
\Stek_{\mix}(\om,g)\sim \Stek(\disk(\mL_S))\sqcup \Stek_N(\hD(\mL_N))\sqcup \Stek_D(\hD(\mL_D)) \sqcup \Stek_{DN}(\qD(\mL_{DN}))
\end{equation}
\end{thm}

Before proving Theorem \ref{asymptotic-dir-Neu}, we state and prove a weaker result.
\begin{lemma}\label{lem for asympt} Theorem~\ref{asymptotic-dir-Neu} holds under the additional assumption that each component of $\pnom$ and of $\pdom$ is geodesic.
\end{lemma}

\begin{proof}
\underline{Case 1.} Suppose $\partial\om=\psom\sqcup \pnom$, so that we have a mixed Steklov-Neumann problem.  Then $\Om$ can be viewed as an orbifold in which each component of $\pnom$ is a reflector.  (See Remark \ref{nota.orb}).   Equation~\eqref{thm:mix asymptotic} in this case is precisely the asymptotic result of \cite{ADGHRS} for orbisurfaces restated in the language of the mixed Steklov-Neumann problem.

\underline{Case 2.} Suppose $\partial\om=\psom\sqcup \pdom$, so that we have a mixed Steklov-Dirichlet problem.  We use the temporary notation $$\pa_*\om:=\pdom \,\,\mbox{and}\,\,\mL_*:=\mL_D.$$
Since $\pa_*\om$ consists of geodesics that are orthogonal to $\psom$ at any common endpoint, we may double $\om$ across $\pa_*\om$ to obtain a Riemannian manifold $\Sigma$ with smooth boundary.  Each component of $\psom$ that is a topological interval becomes a topological circle in $\partial\Sigma$ of twice the length and each topological circle in $\psom$ becomes two topological circles of the same length.  By Equation~\eqref{eq:GPPS} we have 
\begin{equation}\label{eq:gpps2}
\Stek(\Sigma)\sim \Stek(\disk(\mL_S))\sqcup \Stek(\disk(\mL_S))\sqcup \Stek(\disk(2\mL_*)).
\end{equation}

Let $\Stek_N(\om)$ be the mixed Steklov-Neumann spectrum for the same decomposition of $\pa\om$ but with Neumann rather than Dirichlet conditions placed on $\pa_*\om$.  
 Case 2 follows from Equation~\eqref{eq:gpps2} and the following observations:
\begin{itemize}
\item $\Stek(\Sigma)\sim\Stek_D(\om)\sqcup \Stek_N(\om).$
\item $\Stek_N(\om)\sim \Stek(\disk(\mL_S))\sqcup \Stek_N(\hD(\mL_*))$ (by Case 1).
\item $\Stek(\disk(2\mL_*))=\Stek_N(\hD(\mL_*)) \sqcup\,\Stek_D(\hD(\mL_*))$. 
\end{itemize}

\underline{Case 3.} Suppose both $\pdom$ and $\pnom$ are non-trivial.  We may assume that $\om$ satisfies the additional condition that at any point where $\pdom$ and $\pnom$ meet, the interior angle is strictly less than $\pi$.   Indeed, otherwise we may replace $g$ by a conformally equivalent metric $g'$ satisfying this additional condition, where the conformal factor is equal to 1 on a neighborhood of $\pa_S\Om$ as in Remark \ref{add corners}.  This conformal change does not change the mixed spectrum.

The proof of Case 3 is similar to that of Case 2.  We use the temporary notation
$$\mL_1:=\mL_D,\,\,\,\mL_2:=\mL_N,\,\,\,\mL_3:=\mL_{DN}$$

Double $\om$ over $\pnom$ to obtain a Riemannian surface $\Sigma$ with boundary.  Consider the mixed Steklov-Dirichlet problem on $\Sigma$ where the Steklov and Dirichlet parts of the boundary correspond under the doubling to $\psom$ and $\pdom$, respectively.  Then the boundary data for $\Sigma$ can be expressed in terms of the boundary data of $\om$ by:
\begin{equation}\label{bd sigma} \mL_S(\Sigma)=\mL_S\sqcup \mL_S\sqcup 2\mL_2,\,\,\,\,\mL_D(\Sigma)=\mL_1\sqcup \mL_1\sqcup 2 \mL_{3},\,\,\,\,\mL_2(\Sigma)=\mL_{3}(\Sigma)=\emptyset.\end{equation} 

We also consider the Steklov-Dirichlet problem on $\om$ where now we place Dirichlet conditions on \emph{both} $\pdom$ and $\pnom$.  We denote its spectrum by $\Stek_D(\om)$.  Case 3 follows from the following:
\begin{itemize}
\item $\Stek_{D}(\Sigma)=\Stek_{DN}(\om)\sqcup \Stek_{D}(\om).$
\item $\Stek_{D}(\Sigma)\sim \Stek(\disk(\mL_S))\sqcup \Stek(\disk(\mL_S))\sqcup\Stek(\disk(2\mL_2))\sqcup \Stek_D(\hD(\mL_1))\sqcup \Stek_D(\hD(\mL_1))\sqcup \Stek_D(\hD(2\mL_{3}))$ by Case 2 and Equation~\eqref{bd sigma}.
\item $\Stek_{D}(\om)\sim \Stek(\disk(\mL_S))\sqcup \Stek_D(\hD(\mL_1))\sqcup \Stek_D(\hD(\mL_2))\sqcup \Stek_D(\hD(\mL_{3}))$ again by Case 2.
\item $\Stek(\disk(2\mL_2))=\Stek_D(\hD(\mL_2))\sqcup\Stek_N(\hD(\mL_2))$.
\item $\Stek_D(\hD(2\mL_{3}))=\Stek_D(\hD(\mL_{3}))\sqcup\Stek_{DN}(\qD(\mL_3))$ by Equations~\eqref{qD vs hD} and \eqref{hD vs qD+qD}.
\end{itemize}
\end{proof}

\begin{proof}[Proof of Theorem~\ref{asymptotic-dir-Neu}]
By Lemmas~\ref{remove corners} and \ref{lem for asympt}, the theorem holds under the additional hypothesis that each component of $\pdom$ and of $\pnom$ is polygonal.   

Now let $(\om,g)$ and the boundary decomposition $\partial\Om=\psom\sqcup \pnom\sqcup \pdom$ be arbitrary, subject only to the hypotheses of the theorem. We may view $(\om,g)$ as a domain in a complete Riemannian manifold $(M,g)$.  In view of the second hypothesis, we can perturb $\pa \om$ in $M$ to obtain new domains $\om'$ and $\om''$ and corresponding boundary decompositions satisfying the following:
\begin{itemize}
\item $\psom'=\psom=\psom''$ and the intersection of each of $\om$, $\om'$ and $\om''$ with some neighborhood of $\psom$ all coincide.
\item The endpoints (if any) of the components of $\pdom'$ and of $\pdom''$ coincide with those of $\pdom$ and similarly for the Neumann boundary components.
\item Each component of $\pdom'$, $\pnom'$, $\pdom''$, and $\pnom''$ is polygonal.   
\item $\pdom'$, respectively, $\pdom''$, lies in the closure of the exterior, respectively interior, of $\om$.
\item $\pnom'$, respectively, $\pnom''$, lies in the closure of the interior, respectively exterior, of $\om$.
\end{itemize}
Observe that, with these boundary decompositions, $\om$, $\om'$ and $\om''$ all have the same boundary data $(\mL_S,\mL_D,\mL_N,\mL_{DN})$ since we haven't changed the lengths or the types of the Steklov boundary components.    
By domain monotonicity of eigenvalues, we have 
$$\sigma_k^{\mix}(\om',g)\leq \sigma_k^{\mix}(\om,g)\leq\sigma_k^{\mix}(\om'',g).$$
Since $(\om',g)$ and $(\om'',g)$ both satisfy Equation~\eqref{thm:mix asymptotic}, so does $(\om,g)$.
\end{proof}

 \begin{remark}\label{rem: SD and SN asympt} As a special case, consider the mixed Steklov-Neumann and Steklov-Dirichlet problems.  Write $\partial\Om =\psom\sqcup \pst\Om$
 where we will place either Dirichlet or Neumann boundary conditions on $\pst\om$.     We assume this boundary decomposition satisfies the hypotheses of Theorem~\ref{asymptotic-dir-Neu}.
  The theorem together with the model example~\ref{models} gives \[ \sigma_{k+m}^N(\Omega)- \sigma_k^D(\Omega)=O(k^{-\infty})\]
 where $m$ is the number of components of $\psom$ that are topological intervals.      

Specializing further, if each of $\psom$ and $\pa_*\om$ is a topological interval, then we have 
\begin{equation}
\sigma_k^{N}(\Om)L(\pa_S\Om) - \sigma_k^{N}(\hD)L(\hD)  =O(k^{-\infty})
\end{equation}
and 
\begin{equation}
\sigma_k^{D}(\Om)L(\pa_S\Om) - \sigma_k^{D}(\hD)L(\hD)  =O(k^{-\infty}).
\end{equation}

If in addition, $\om$ is assumed to be a simply-connected plane domain and $\psom$ is assumed to be a straight line segment of length $\ell$, Levitin, Parnovski, Polterovich, and Sher \cite[Propositions 1.3 and 1.8]{LPPS}, showed that there exists a constant $c>0$ such that
\[\sigma_{k}^N(\Omega)= \sigma_k^N(\hD(\ell))+O(e^{-ck}), \quad\sigma_{k}^D(\Omega)= \sigma_k^D(\hD(\ell))+O(e^{-ck}).\]   

They also considered the situation in which the angle between $\pa_S\Om$ and $\pa_*\Om$ is not $\pi/2$, obtaining a two-term asymptotic that verifies a conjecture in \cite{FK}. Note that it is unlikely that we can get exponential decay in our general setting in which $\psom$ is assumed only to be smooth, not even analytic. 
 \end{remark}
 
 We give several applications of Theorem~\ref{asymptotic-dir-Neu}.

\begin{cor}\label{cor: inverse} Let $(\Om,g)$ be a connected, compact Riemannian surface with a boundary decomposition $\partial\Om=\psom\sqcup \pa_*\om$.  Consider the mixed Steklov-Dirichlet, respectively mixed Steklov-Neumann, problem on $\om$ with $\pa_*\om$ playing the role of $\pdom$, respectively $\pnom$.    Assume the hypotheses of Theorem~\ref{asymptotic-dir-Neu} hold.   Let $n$ and $m$ be the number of components of $\psom$ that are topological circles and topological intervals, respectively and let $\mL_S$ and $\mL_*$ be the multisets of their lengths.   Then the following are spectral invariants of each of $\Stek_D(\om,g)$ and $\Stek_N(\om,g)$:

  \begin{itemize}
 \item $m$ and $n$;
 \item $\mL_S\sqcup\mL_S\sqcup 2\mL_*$.
 \end{itemize}

\end{cor}

\begin{remark}\label{rem: 5.7} Let $\mathcal{L}$ be the class of all pairs of finite multisets $(\mL_S,\mL_*)$ of positive real numbers.  Given $(\mL_S,\mL_*)\in \mathcal{L}$ such that $\mL_*$ has a repeated entry, define an ``entry-exchange'' as follows:  Choose $\ell\in \mL_S$ and choose $\ell^*\in \mL_* $ such that $\ell^*$ has multiplicity at least two.    In $\mL_S$, replace $\ell$ by $2\ell^*$ and, in $\mL_*$, replace two copies of $\ell^*$ by two copies of $\frac{\ell}{2}$ to get a new pair $(\mL'_S,\mL'_*)\in \mathcal{L}$.    The operation ``entry-exchange'' generates an equivalence relation on $\mathcal{L}$.  (In particular, if $\mL_*$ has no repeated entries, then the equivalence class $[(\mL_S,\mL_*)]$ consists only of $(\mL_S,\mL_*)$.)  Corollary~\ref{cor: inverse} says that for Riemannian surfaces with boundary decompositions satisfying the stated hypotheses, each of the Steklov-Dirichlet spectrum and the Steklov-Neumann spectrum determines the boundary data $(\mL_S,\mL_*)$ up to ``entry-exchange'' equivalence.
   \end{remark}

\begin{cor} We use the notation and hypotheses of Corollary~\ref{cor: inverse}.
\begin{enumerate}
\item[(i)] Assume $\Omega$ is simply-connected and $m=1$ (and thus $n=0$).  Then for every $\epsilon>0$ there exists $k_{\Omega}\in\mathbb{N}$ such that for any $k\ge k_\Om$ we have 
 \[\max\{\sigma_{k}^N(\Omega), \sigma_{k-1}^D(\Om)\}L(\pa_S\Om)< (2k-1)\pi+\epsilon.\]
  \item[(ii)] More generally, if $\Omega$ has genus zero and if all components of $\pa\om$ except for one lie entirely in $\psom$,  then there exist $k_\Om>0$ such that for $k>k_\Om$ 
  \[\sigma_{k}^N(\Omega)L(\pa_S\Om)< 4(2k-1)\pi+\epsilon,\qquad \sigma_{k}^D(\Om)L(\pa_S\Om)< 4(2(k+m)-1)\pi+\epsilon.\] 
  \end{enumerate}
\end{cor}
\begin{proof}

By Remark~\ref{rem: SD and SN asympt}, we can choose $k_{\Omega}$ such that for every $k\ge k_\Om$ we have  
\begin{equation}\label{eq:kom} |\sigma_{k}^N(\Om)-\sigma_{k-m}^D(\Om)|< \frac{\epsilon}{L(\psom)}.\end{equation} 

Statement (i) follows from Equation~\eqref{eq:kom} and Theorem~\ref{general minpi}, since $\Sigma_{2k-1}(0,1)\leq 2(2k-1)\pi$.

For (ii), Theorem~\ref{minpiwithholes} yields
$$\min\{\sigma_{k}^N(\Omega), \sigma_{k-m}^D(\Om)\}L(\pa_S\Om)\leq \min\{\sigma_{k}^N(\Omega), \sigma_{k-1}^D(\Om)\}L(\pa_S\Om)< 4(2k-1)\pi.$$
We can then apply Equation~\eqref{eq:kom} to complete the proof.

\end{proof}

Our final corollary is an asymptotic bound on the  multiplicity of eigenvalues for all mixed Steklov problems satisfying the hypotheses of Theorem~\ref{asymptotic-dir-Neu}.   Our results extend the asymptotic multiplicity bounds in  \cite[Corollary 1.6]{GPPS.14} (see also \cite[Propostion 1.5]{KKP.14}).  They showed for the pure Steklov problem that the multiplicity of $\sigma_k(\Om)$  is bounded above asymptotically by $2l$, where $l$ is the number of boundary components of $\pa\Om$.

\begin{cor}\label{asymtotic mult} Let $(\Om,g)$ and the decomposition $\pa\om=\psom\sqcup\pdom\sqcup\pnom$ satisfy the assumptions of Theorem \ref{asymptotic-dir-Neu}.  Let $n$ and $m$ be the number of components of $\psom$ that are topological circles and topological intervals, respectively.  Then:
\begin{itemize} \item[(i)]  there exists $k_\Om>0$, such that the multiplicity of $\sigma_k^{\mix}(\Om)$ is at most $2n+m$  for all $k \ge k_\Om $. 
 \item[(ii)] If the lengths $\ell_1,\ell_2,\dots, \ell_{n+m}$ of the components of $\psom$ are algebraic numbers and none is a rational multiple of any other, then there exists $k_\om>0$, such that the multiplicity of $\sigma_k^{\mix}(\Om)$ is at most $2$ for all $k>k_\Om$.  Moreover, if $n=0$, then $\sigma_k^{\mix}(\Om)$ is simple for all $k>k_\Om$. \end{itemize}
\end{cor}

\begin{proof}[Proof of Corollary \ref{asymtotic mult}]

   (i) is an immediate consequence of Theorem \ref{asymptotic-dir-Neu} and the fact that the eigenvalues in $\Stek(\disk(\mL_S))\sqcup \Stek_N(\hD(\mL_N))\sqcup \Stek_D(\hD(\mL_D)) \sqcup \Stek_{DN}(\qD(\mL_{DN}))$ have multiplicity at most $2n+m$.  (See Example~\ref{models}.)
 
  (ii) We use the fact that if $x,y$ are algebraic numbers and not rational multiples of each other, then there exist positive constants $c(x,y)$ and $d(x,y)$  depending only on $x$ and $y$ such that 
   \begin{equation}\label{OK}|px-q y|\ge\frac{c(x,y)}{(pq)^{d(x,y)}},\end{equation} for all $p,q\in\mathbb{N}$.  In particular, $d(x,y)$ is an integer depending only on the degree of the minimal polynomial of $x$ and $y$.  We thank Oleksiy Klurman for pointing out the inequality above to us.
   
    Now, let $x_i:=\frac{1}{2\ell_i}$, $1\le i\le n+m$.   Each pair $x_i,x_j$, $i\neq j$, satisfies inequality \eqref{OK} for some $c(x_i,x_j), d(x_i,x_j)>0$. Let $$c=\max_{i,j}\{c(x_i,x_j)\},\qquad d=2\max_{i,j}\{d(x_i,x_j)\}.$$ 
    Hence, for any $p,q\in \mathbb{N}$ and $i\ne j$, we have $|\frac{p}{\ell_i}-\frac{q}{\ell_j}|\ge\frac{c}{\max\{p,q\}^{d}}$.  Item (ii) now follows from Theorem~\ref{asymptotic-dir-Neu} and Example~\ref{models}.
     \end{proof}

\bibliographystyle{plain}
\bibliography{refpp2.bib}

\end{document}